\documentclass[aop,MSNbibl,nameyear,dvips]{arximspdf}
\usepackage{mathrsfs}
\usepackage{graphicx}

\doi{10.1214/11-AOP658}
\volume{40}
\issue{4}
\pubyear{2012}
\firstpage{1636}
\lastpage{1674}

\makeatletter

\newtheorem{theorem}{Theorem}
\newtheorem{corollary}{Corollary}
\newtheorem{proposition}{Proposition}
\newtheorem{lemma}{Lemma}

\newproclaim{defi}{Definition}
\newproclaim{remark}{Remark}

\makeatother

\begin{document}
\begin{frontmatter}

\title{The convex minorant of a L\'{e}vy process}
\runtitle{The convex minorant of a L\'{e}vy process}
\pdftitle{The convex minorant of a Levy process}

\begin{aug}
\author{\fnms{Jim} \snm{Pitman}\ead[label=e1]{pitman@stat.Berkeley.edu}\thanksref{t1}}
\and
\author{\fnms{Ger\'{o}nimo} \snm{Uribe Bravo}\corref{}\ead[label=e2]{geronimo@matem.unam.mx}\thanksref{t2}}

\pdfauthor{Jim Pitman and Geronimo Uribe Bravo}
\runauthor{J. Pitman and G. Uribe Bravo}
\affiliation{University of California, Berkeley}
\address{Department of Statistics\\
University of California, Berkeley\\
367 Evans Hall\\
Berkeley, California 94720-3860\\
USA\\
\printead{e1}\\
\phantom{E-mail: }\printead*{e2}} 
\thankstext{t1}{Supported in part by NSF Grant DMS-08-06118.}
\thankstext{t2}{Supported by a~postdoctoral fellowship from UC
MexUS---CoNaCyt and NSF Grant DMS-08-06118.}
\end{aug}

\received{\smonth{11} \syear{2010}}
\revised{\smonth{2} \syear{2011}}

%
\begin{abstract}
We offer a~unified approach to the theory of convex minorants of L\'
{e}vy processes with continuous distributions. New results include
simple explicit constructions of the convex minorant of a~L\'{e}vy
process on both finite and infinite time intervals, and of a~Poisson
point process of excursions above the convex minorant up to an
independent exponential time. The Poisson--Dirichlet distribution of parameter 1 is shown to be
the universal law of ranked lengths of
excursions of a~L\'{e}vy process with continuous distributions above its
convex minorant on the interval
$[0,1]$.
\end{abstract}

%
\setattribute{keyword}{AMS}{MSC2010 subject classification.}
\pdfkeywords{60G51, Levy processes, convex minorant, uniform stick-breaking, fluctuation theory, Vervaat transformation}

\end{frontmatter}
%

\section{Introduction}
\label{IntroductionSection}

We present simple explicit constructions of the convex minorant of a~L\'
{e}vy process with continuous distributions on both finite and infinite
time intervals, and of a~Poisson point process of excursions of the L\'
{e}vy process above its convex minorant. These constructions bridge
a~number of gaps in the literature by relating combinatorial approaches
to fluctuation theory of random walks related to the cycle structure of
random permutations, dating back to the 1950s [cf.~\citeauthor{MR0039939}
(\citeyear{MR0039939,MR0060173,MR0058893,MR0068154});
\citet{MR0079851}], some features of
which were extended to interval partitions associated with the convex
minorant of Brownian motion and Brownian bridge by
\citeauthor{MR1978665} (\citeyear{MR1978665,MR1861441}) and \citet{Balabdaoui2009fk},
and results previously obtained for the convex minorants of Brownian motion
by \citet{MR714964} and \citet{MR733673}, and for L\'{e}vy processes by
\citet{MR1739699} and \citet{MR1747095}. In particular,
we gain access
to the excursions above the convex minorant, which were previously
treated only in the Brownian case by \citet{MR714964} and
\citet{MR733673}.

Our work is part of a~larger initiative to understand the convex
minorant of processes with exchangeable\vadjust{\goodbreak} increments. The case of
discrete time is handled in \citet{2010AbramsonPitman}, while Brownian
motion is given a~more detailed study in \citet
{2010PitmanRoss}. Our
joint findings are summarized in \citet{Abramson2011fk}.

\subsection{Statement of results}
Let $X$ be a~L\'{e}vy process. The following hypothesis is used
throughout the paper:
\begin{enumerate}[(CD)]
\item[(CD)]\label{CD}For all $t>0$, $X_t$ has a~continuous
distribution, meaning that for each $x\in\mathbb{R}$, $\mathbb{P} (
X_t=x )=0$.
\end{enumerate}
It is sufficient to assume that $X_t$ has a~continuous distribution for
some $t>0$. Equivalently [\citeauthor{MR1739520} (\citeyear{MR1739520}), Theorem 27.4, page
175] 
$X$ is not a~compound Poisson process with drift.

The \textit{convex minorant} of a~function $f$ on an interval
$[0,t]$ or
$[0,\infty)$ is the greatest convex function $c$ satisfying $c\leq f$.
We shall only consider functions $f$ which are c\`{a}dl\`{a}g, meaning that
$\lim_{h\to0+}f ( t+h) =f ( t) $ and that $\lim_{h\to0-}f ( t-h) $
exists; the latter limit will be denoted $f ( t-) $. 

First properties of the convex minorant of a~L\'{e}vy process,
established in Section~\ref{BasicPropertiesSection} and which partially
overlap with the Markovian study of convex minorants in \citet
{Lachieze-Rey2009fk}, are:
\begin{proposition}
\label{BasicPropertiesOfConvexMinorant}
Let $X$ be a~L\'{e}vy process which satisfies \textup{(CD)} and $C$ the
convex minorant of $X$ on $[0,t]$. The following conditions hold almost surely:
\begin{longlist}[1.]
\item[1.]\label{ComplementIsLightProperty} The open set $\mathscr{O}=\{
s\in (0,t)\dvtx C_s<X_s\wedge X_{s-}\} $ has Lebesgue measure $t$.
\item[2.]\label{OneJumpPerExcursionProperty} For every component interval
$( g,d) $ of $\mathscr{O}$, the jumps that $X$ might have at~$g$ and $d$
have the same sign.
When $X$ has unbounded variation on finite intervals, both jumps are zero.
\item[3.]\label{SimplicityOfSlopesProperty}If $(g_1,d_1)$ and $( g_2,d_2)
$ are different component intervals of $\mathscr{O}$, then their
slopes differ
\[
\frac{C_{d_1}-C_{g_1}}{d_1-g_1}\neq\frac{C_{d_2}-C_{g_2}}{d_2-g_2}.
\]
\end{longlist}
\end{proposition}

Let $\mathscr{I}$ be the set of connected components of $\mathscr
{O}$; we shall
also call them \textit{excursion intervals}. Associated with each
excursion interval $( g,d) $ are the \textit{vertices} $g$ and $d$,
the \textit{length} $d-g$, the \textit{increment} $C_d-C_g$ and
the \textit{slope} $( C_d-C_g) /( d-g) $.

Our main result is a~simple description of the lengths and increments
of the excursion intervals of the convex minorant.
Indeed, we will consider a~random ordering of them which uncovers
a~remarkable probabilistic structure.

\begin{theorem}
\label{PointProcessDescriptionForDeterministicAndFiniteHorizon}
Let $( U_i) $ be a~sequence of uniform random variables on $(0,t)$
independent of the L\'{e}vy process $X$ which satisfies \textup{(CD)}.
Let $( g_1,d_1) $, $( g_2,d_2) ,\ldots$ be the sequence of
distinct excursion intervals which are successively discovered by the
sequence $( U_i) $. Consider another i.i.d. sequence $( V_i) $
of uniform random variables on $(0,1)$ independent of $X$, and
construct the associated \textup{uniform stick-breaking process} $L$ by
\[
L_1=tV_1 \quad\mbox{and}\quad \mbox{for }i\geq1\qquad L_{i+1}=V_{i+1}( t-S_i),
\]
where
\[
S_0=0  \quad\mbox{and}\quad \mbox{for }i\geq1\qquad S_{i}=L_1+\cdots+L_i.
\]
Under hypothesis \textup{(CD)}, the following equality in
distribution holds:
\[
\bigl( ( d_i-g_i,C_{d_i}-C_{g_i}) , i\geq1\bigr) \stackrel{d}{=}\bigl( ( L_i,
X_{S_i}-X_{S_{i-1}}) ,i\geq1\bigr) .
\]
\end{theorem}
%


The \textit{Poisson--Dirichlet distribution of parameter one} is the law
of the decreasing rearrangement of the sequence $L$ when $t=1$. Theorem
\ref{PointProcessDescriptionForDeterministicAndFiniteHorizon} implies
that the Poisson--Dirichlet distribution of parameter 1 is the
universal distribution of the ranked lengths of excursions intervals of
the convex minorant of a~L\'{e}vy process with continuous distributions
on $[0,1]$.
%
%
What differs between each L\'{e}vy process is the
distribution of the order in which these lengths appear, that is, the
law of the composition of of $[0,1]$ induced by the lengths of
excursion intervals when they are taken in order of appearance.
Using Theorem~\ref{PointProcessDescriptionForDeterministicAndFiniteHorizon} we can form a~composition of $[0,1]$ with that law in the following way. For each
pair $( L_i,X_{S_i}-X_{S_{i-1}}) $ we generate a~slope by dividing
the second coordinate, the increment, by the first, the length, and
then create a~composition of $[0,1]$ by arranging the sequence $L$ in
order of increasing associated slope.

Note that the second sequence of Theorem~\ref{PointProcessDescriptionForDeterministicAndFiniteHorizon} can also be
constructed as follows: given a~uniform stick-breaking process $L$,
create a~sequence $Y_i$ of random variables which are conditionally
independent given $L$ and such that the law of $Y_i$ given $L$ is that
of $X_{L_i}$ ($X$ independent of $L$). Then
\[
\bigl( ( L_i,Y_i) \dvtx i\geq1\bigr) \stackrel{d}{=}\bigl( ( L_i, X_{S_i}-X_{S_{i-1}})
,i\geq1\bigr) .
\]

Theorem~\ref{PointProcessDescriptionForDeterministicAndFiniteHorizon}
provides a~way to perform explicit computations. For example, the
intensity measure $\nu_t$ of 
the point process $\Xi_t$ with atoms at
\[
\{ ( d-g,C_d-C_g) \dvtx ( g,d) \mbox{ is an excursion interval}\}
\]
is given by
%
\begin{eqnarray*}
\nu_t ( A)
&=&\mathbb{E} \biggl( \sum_i {\mathbf{1}}_{( d_i-g_i,C_{d_i}-C_{g_i}) \in
A} \biggr)
\\
&=&\mathbb{E} \biggl( \sum_i {\mathbf{1}}_{( L_i,X_{d_i}-X_{g_i}) \in A} \biggr)
=\int_0^t \int{\mathbf{1}}_{A} ( l,x) \frac{1}{l} \mathbb{P} (
X_l\in dx ) \,dl.
\end{eqnarray*}
[This follows conditioning on $L$ and then using the intensity measure
of $L$ obtained by size-biased sampling; cf. formula (6) in
\citet{MR1434129}.]


We now apply Theorem~\ref{PointProcessDescriptionForDeterministicAndFiniteHorizon} to fully
describe the convex minorant of the Cauchy process as first done in
\citet{MR1747095}.
Let $X$ be a~Cauchy process characterized by
\[
F ( x) : =\mathbb{P} ( X_1\leq x )=1/2+\arctan ( x) /\pi.
\]
Let $C$ be the convex minorant of $X$ on $[0,1]$ and $D$ its right-hand
derivative,
\[
D_t=\lim_{h\to0+}\frac{C_{t+h}-C_t}{h}.
\]
Consider
\[
I_x=\inf\{ t\geq0\dvtx  D_t>x\}\qquad \mbox{for $x\in\mathbb{R}$}.
\]
Note that $\mathbb{P} ( X_t<xt )=F ( x) $ and that therefore, in the
setting of Theorem~\ref{PointProcessDescriptionForDeterministicAndFiniteHorizon}, the slopes
$(C_{d_i}-C_{g_i})/(d_i-g_i)$ are independent of the lengths $d_i-g_i$.
Also, let $T$ be a~Gamma subordinator such that
\[
\mathbb{E} ( e^{-q T_t} )=\biggl( \frac{1}{1+q}\biggr) ^t.
\]

\begin{corollary}
\label{CauchyCorollary}
\label{CauchyCharacterization}\textup{1.} The symmetric Cauchy process is
characterized by the independence of lengths and slopes of excursions
intervals on $[0,1]$.
\begin{longlist}[2.]
\item[2.]\label{BertoinsTheorem} $( I_x,x\in\mathbb{R}) $ and $( T_{F (
x) }/T_1,x\in\mathbb{R}) $ have the same law.
\end{longlist}
\end{corollary}

Item \hyperref[BertoinsTheorem]{2} is due to \citet{MR1747095}, who
used a~technique allowing only the study of the convex minorant of a~Cauchy
process on $[0,1]$.



Integrating Theorem~\ref{PointProcessDescriptionForDeterministicAndFiniteHorizon}, we obtain a~description of the convex minorant considered on the random interval
$[0,T_\theta]$ where $T_\theta$ is a~exponential random variable of
parameter $\theta$ independent of $X$.
\begin{corollary}
\label{PointProcessDescriptionForIndependentExponentialHorizon}
Let $T$ be exponential of parameter $\theta$ and independent of the L\'
{e}vy process $X$ which satisfies \textup{(CD)}. Let $\Xi_{T}$ be the
point process with atoms at lengths and increments of excursion
intervals of the convex minorant of $X$ on $[0,T]$. Then $\Xi_{T}$ is a~Poisson point process with intensity
\[
\mu_\theta ( dt,dx) =e^{-\theta t}\frac{dt}{t}\mathbb{P} ( X_t\in
dx ).
\]
\end{corollary}

By conditioning on $T$ (which essentially reduces to inverting Laplace
transforms and underlies the analysis of the relationship between the
Gamma subordinator and the Poisson--Dirichlet distribution), we see
that Theorem~\ref{PointProcessDescriptionForDeterministicAndFiniteHorizon} can be
deduced from Corollary~\ref{PointProcessDescriptionForIndependentExponentialHorizon}. The latter
can be deduced from the analysis of the independence of pre- and
post-minimum processes of a~L\'{e}vy process run until\vadjust{\goodbreak} an independent
exponential time found in \citet{MR588409}. 
These relationships are discussed in Section~\ref{SplittingAtTheMinimumSection}, where we also explain the results on
fluctuation theory for L\'{e}vy processes which are found in the
literature and which can be deduced from our analysis of the convex
minorant.\looseness=-1 

From Theorem~\ref{PointProcessDescriptionForDeterministicAndFiniteHorizon} we can also
derive the behavior of the convex minorant of $X$ on $[0,\infty)$ as
described for a~Brownian motion by \citet{MR714964} %
and \citet{MR733673} %
and for a~L\'{e}vy process by \citet{MR1739699}. Let $\Xi
_\infty$ be
the point process of lengths of excursion interval and increments of
the convex minorant on $[0,\infty)$.
\begin{corollary}
\label{PointProcessDescriptionForInfiniteHorizon}
The quantity $l=\liminf_{t\to\infty}X_t/t$ belongs to $(-\infty
,\infty
]$ and is almost surely constant if and only if the convex minorant of
$X$ on $[0,\infty)$ is almost surely finite.
In this case, under \textup{(CD)}, $\Xi_\infty$ is a~Poisson point
process with intensity
\[
\mu_\infty ( dt,dx) =\frac{{\mathbf{1}}_{x<l t}}{t} \mathbb{P} (
X_t\in dx ) \,dt.
\]
\end{corollary}

Recall, for example, \citeauthor{110460001} [(\citeyear{110460001}), Example 7.2], the
strong law of
large numbers for L\'{e}vy processes, which says that if the
expectation of $X_1$ is defined, then
\[
\lim_{t\to\infty} \frac{X_t}{t}=\mathbb{E} ( X_1 ) \qquad\mbox{almost surely}.
\]
Hence, if $\mathbb{E} ( X_1^- )<\infty$, we can apply the second part of
Corollary~\ref{PointProcessDescriptionForInfiniteHorizon} with
$l=\mathbb{E} ( X_1 )$. In the remaining case when $\mathbb{E} (
X_1^- )=\mathbb{E} ( X_1^+ )=\infty$,
let $\nu$ be the L\'{e}vy measure of $X$ and $\nu_+$ its right-tail
given by
\[
\nu_+ ( y) =\nu ( ( y,\infty) ) .
\]
\citet{MR0336806} provides the necessary and sufficient for
$-\infty<l$,
which implies that, actually, $l=\infty$,
\[
\int_{(-\infty,0)}\frac{|y|}{\nu_+ ( |y|) } \nu ( dy) <\infty
\]
(see also \citeauthor{MR2320889} (\citeyear{MR2320889}), page 39, for a~proof).

While it seems natural to first study the convex minorant of a~L\'{e}vy
process on $[0,\infty)$, as was the approach of previous authors, the
description of the convex minorant with
infinite horizon is less complete, as it is necessarily restricted to
slopes $a < l$.

As another application, we can use the stick-breaking representation of
Theorem~\ref{PointProcessDescriptionForDeterministicAndFiniteHorizon}
to study the absolute continuity of the location and the value of the
minimum of the L\'{e}vy process on $[0,1]$. Let
\[
\underline X_t=\min_{s\leq t}X_s \quad\mbox{and}\quad\rho_t\mbox{ be such
that } X_{\rho_t}\wedge X_{\rho_t-}=\underline X_t.
\]
(Recall that under {(CD)}, the minimum of a~L\'{e}vy process on
$[0,t]$ is attained at an almost surely unique place $\rho_t$, as
deduced from Theorem~\ref{PointProcessDescriptionForDeterministicAndFiniteHorizon} since
$\mathbb{P} ( X_t=0 )=0$.)\vadjust{\goodbreak}
\begin{theorem}
\label{AbsoluteContinuityTheorem}
Let $X$ be a~L\'{e}vy process such that $0$ is regular for both
half-lines $(-\infty,0)$ and $(0,\infty)$. Then:
\begin{longlist}[1.]
\item[1.] The distribution of $\rho_1$ is equivalent to Lebesgue measure
on $[0,1]$.
\item[2.] If $X_t$ has an absolutely continuous distribution for each
$t>0$, then the distribution of $( \rho_1,\underline X_1) $ is
equivalent to Lebesgue measure on $(0,1]\times(0,\infty)$.
\item[3.] If $X_t$ has an absolutely continuous distribution for each
$t>0$, then the distribution of $( \underline X_1,X_1-\underline X_1)
$ is equivalent to Lebesgue measure on $(-\infty,0)\times(0,\infty)$.
\end{longlist}
\end{theorem}

\citet{Chaumont2010fk} also analyzes absolute continuity properties
for the supremum of a~L\'{e}vy process on a~fixed interval using
excursion theory for the reflected L\'{e}vy process. The densities
provided by Theorem~\ref{AbsoluteContinuityTheorem} (more importantly,
the fact that they are almost surely positive) provide one way to
construct bridges of the L\'{e}vy process $X$ conditioned to stay
positive. With these bridges, we can prove a~generalization of
Vervaat's theorem relating the Brownian bridge and the normalized
Brownian excursion [\citeauthor{MR515820} (\citeyear{MR515820}), Theorem 1] to a~fairly
general class
of L\'{e}vy processes. Details are provided in \citet
{UribeBravoVervaatExtended}.

Our next results will only consider convex minorants on a~fixed
interval, which we take to be $[0,1]$.

Theorem~\ref{PointProcessDescriptionForDeterministicAndFiniteHorizon}
gives a~construction of the convex minorant by means of sampling the L\'
{e}vy process at the random, but independent, times of a~uniform
stick-breaking process. Our second proof of it, which does not rely on
fluctuation theory and gives insight into the excursions of $X$ above
its convex minorant, depends on the use of the following path
transformation. Let $u$ be an element of the excursion set $\mathscr{O}$,
and let $( g,d) $ be the excursion interval which contains $u$. We
then define a~new stochastic process $X^u=( X^u_t) _{t\leq1}$ by
%
%
\begin{equation}
X^u_t=
\cases{
X_{u+t}-X_u,&\quad$0\leq t< d-u$,\vspace*{2pt}\cr
C_d-C_g+X_{g+t-(d-u)}-X_{u},&\quad$d-u\leq t\leq d-g$,\vspace*{2pt}\cr
C_d-C_g+X_{t-(d-g)},&\quad$d-g\leq t<d$,\vspace*{2pt}\cr
X_{t},&\quad$\mbox{$d\leq t\leq1$}$.}
\label{PathTransformationForLevyProcess}
\end{equation}
\begin{figure}

\includegraphics{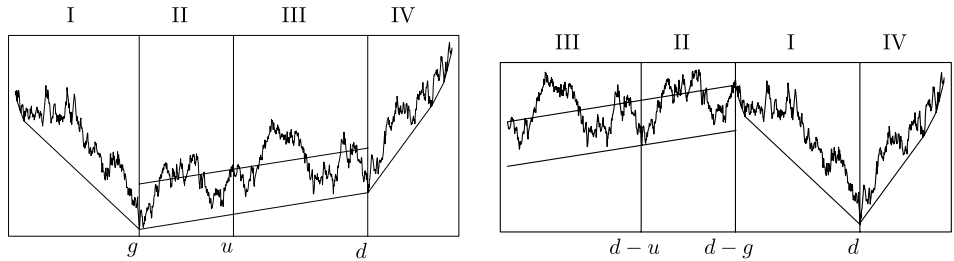}

\caption{Visualization of the path transformation $X\mapsto X^u$
applied to a~Brownian motion seen from its convex minorant.}
\label{PathTransformation}
\end{figure}


\begin{figure}

\includegraphics{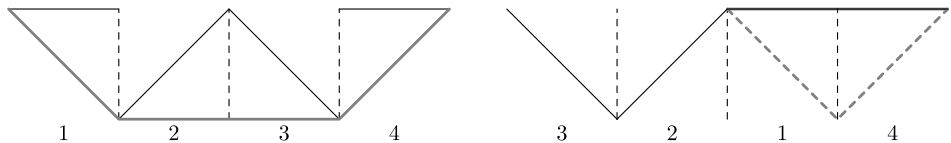}

\caption{Visualization of the path transformation $X\mapsto X^u$
applied to a~c\`{a}dl\`{a}g path not satisfying property 2 of Proposition \protect\ref{BasicPropertiesOfConvexMinorant}.}
\label{PathTransformation2}
\end{figure}



The idea for such a~definition is that the graph of the convex minorant
of~$X^u$ on $[d-g,1]$ can be obtained from the graph of $C$ by removing
$( g,d) $ and closing up the gap adjusting for continuity, while on
$[0,d-g]$, $X^u$ goes from $0$ to $C_d-C_g$. (Property \hyperref[OneJumpPerExcursionProperty]{2} of Proposition~\ref{BasicPropertiesOfConvexMinorant} is essential for the transformation
to work like this; see Figure~\ref{PathTransformation2}.) A~schematic
picture of the path transformation is found in Figure~\ref{PathTransformation} for a~typical Brownian motion path.

Theorem~\ref{PointProcessDescriptionForDeterministicAndFiniteHorizon}
then follows from the following invariance result. Indeed, by applying
the following path transformation recursively, we can obtain a~size-biased sample of the excursion intervals. In particular, the
excursion interval containing an independent uniform variable has a~uniform length, which begins to explain the stick-breaking process
appearing in Theorem~\ref{PointProcessDescriptionForDeterministicAndFiniteHorizon}.
\begin{theorem}
\label{InvarianceUnderPathTransformation}
If $U$ is a~uniform on $(0,1)$ and independent of $X$ and hypothesis
\textup{(CD)} holds, the pairs $( U,X) $ and $( d-g,X^U) $ have
the same law.
\end{theorem}

Proof of Theorem~\ref{InvarianceUnderPathTransformation} will be based
on the analogous random walk result proved by \citet
{2010AbramsonPitman} as well as analysis on Skorohod space.
Abraham and Pitman's discrete time result is an exact invariance
property for a~similar transformation applied to the polygonal
approximation $X^n$ of $X$ given by $X^n_t=X_{[nt]/n}( \lceil nt \rceil
/n-t) +X_{\lceil nt \rceil/n}( [nt]/n) $; that this approximation
does not converge in Skorohod space to $X$ makes the passage to the
limit technical, although it simplifies considerably for L\'{e}vy
processes with unbounded variation, and particularly so for L\'{e}vy
processes with continuous sample paths. The discrete time result is
combinatorial in nature and related to permutations of the increments.
Indeed, the discrete time result is based on the fact that for a~random
walk $S$ with continuous jump distribution, the probability that $S$
lies strictly above the line from $(0,0)$ to $( n,S_n) $ on $\{
0,\ldots, n\} $ is known to be $1/n$, and conditionally on this event,
the law of $S$ can be related to a~Vervaat-type transform of $S$.
Hence, it is not only possible to verify by combinatorial reasoning
that the faces of the convex minorant have the same law as the cycle
lengths of a~uniform random permutation when both are placed in
decreasing order, but also to characterize the path fragments on top of
each excursion interval.

Theorem~\ref{InvarianceUnderPathTransformation} actually gives a~much
stronger result than Theorem~\ref{PointProcessDescriptionForDeterministicAndFiniteHorizon} since it
grants access to the behavior of $X$ between vertex points of the
convex minorant. To see this, consider the Vervaat transformation: for
each \mbox{$t>0$} and each c\`{a}dl\`{a}g\ function~$f$, let $\rho_t=\rho_t
( f) $
be the location of the last mini\-mum~$\underline f ( t) $ of $f$ on
$[0,t]$ and define
\[
V_tf ( s) =f ( \rho_t+s\operatorname{mod} t) -\underline{f} ( t) \qquad\mbox{for
$s \in[0,t]$}.
\]
This path transformation was introduced in \citet{MR515820} for the
Brownian bridge; its connection to L\'{e}vy processes was further
studied for stable L\'{e}vy processes by \citet{MR1465814},
for spectrally positive L\'{e}vy processes in \citet{MR1844511},
and more general L\'{e}vy processes by \citet{MR2139029}.

For each excursion interval $( g,d) $ of $\mathscr{O}$, associate an
\textit{excursion} $e^{( g,d) }$ given by
\[
e^{( g,d) } ( s) =X_{g+s}-C_{g+s}\qquad\mbox{for $s\in[0,d-g]$;}
\]
note that $e^{( g,d) } ( 0) $ is positive if $X_{g}>C_g$.
Finally, recalling the setting of Theorem~\ref{PointProcessDescriptionForDeterministicAndFiniteHorizon}, let $K^i$ be
Knight's bridge,
\[
K^i_s=X_{( S_{i-1}+t) }-X_{S_{i-1}}-s\frac
{X_{S_i}-X_{S_{i-1}}}{L_i},\qquad \mbox{$s\in[0,L_i]$}
\]
[the name is proposed because of remarkable universality theorems
proven for $K^i$ in \citet{MR1417982}].

%
\begin{theorem}
\label{PointProcessOfExcursions}
The following equality in distribution holds under \textup{(CD)}:
%
%
\begin{equation}
\label{EqualityInLawOfPointProcessOfExcursions}
\qquad\bigl( \bigl( d_i-g_i,C_{d_i}-C_{g_i},e^{( g_i,d_i) }\bigr) ,i\geq1\bigr)
\stackrel{d}{=}
\bigl( \bigl( L_i,X_{S_i}-X_{S_{i-1}},V_{L_i} ( K^i) \bigr) ,i\geq1\bigr) .
\end{equation}
\end{theorem}

Note that the increment $C_d-C_g$ cannot be obtained from the path
fragment $e^{( g,d) }$ when $X$ jumps at $g$ or $d$. This does not
happen if $X$ has unbounded variation, thanks to Proposition~\ref{BasicPropertiesOfConvexMinorant}.

The same remark of Theorem~\ref{PointProcessDescriptionForDeterministicAndFiniteHorizon} holds,
namely, the intensity measure of the right-hand side of (\ref{EqualityInLawOfPointProcessOfExcursions}), seen as a~point process,
admits the expression
\[
\mathbb{E} \biggl( \sum_i{\mathbf{1}}_{( L_i,X_{S_i}-X_{S_{i-1}},V_{L_i} (
K^i) ) \in A} \biggr)=\int_0^1 \int\int{\mathbf{1}}_{A}{( l,x,e) }\frac
{1}{l} \kappa _l ( dx,de) \,dl
\]
in terms of the law of $X$, where the measure $\kappa_l$ is the joint
law of $X_l$ and the Vervaat transform $V_l$ of $( X_{t}-tX_l/l,t\in
[0,l]) $. The measure $\kappa_l$ is related to L\'
{e}vy processes conditioned to stay positive [introduced in generality
in \citet{MR2164035}] in \citet{UribeBravoVervaatExtended}.

This document is organized as follows: we first study the basic
properties of the convex minorant of a~L\'{e}vy process of Proposition
\ref{BasicPropertiesOfConvexMinorant} in Section~\ref{BasicPropertiesSection}. Then, examples of the qualitative behaviors
of the convex minorants are given in Section~\ref{ExamplesSection}.
Next, we turn to the description of the process of lengths and slopes
of excursion intervals up to an independent exponential time in Section
\ref{SplittingAtTheMinimumSection}, where we also discuss how this
implies the description of the convex minorant to a~deterministic and
finite time and on an infinite horizon. Section~\ref{SplittingAtTheMinimumSection} also explains the relationship between
this work and the literature on fluctuation theory for L\'{e}vy
processes. Section~\ref{ApplicationsSection} is devoted to the absolute
continuity of the location and time of the minimum of a~L\'{e}vy
process with a~proof of Theorem~\ref{AbsoluteContinuityTheorem}.
Finally, we pass to the invariance of the path transformation (\ref{PathTransformation}) for
L\'{e}vy processes stated as Theorem~\ref{InvarianceUnderPathTransformation}, in Section~\ref{InvariantPathTransformationSection} and to the description of the
excursions above the convex minorant implied by Theorem~\ref{PointProcessOfExcursions} in Section~\ref{ExcursionSection}.

\section{Basic properties of the convex minorant on a~finite interval}
\label{BasicPropertiesSection}
In this section we will prove Proposition~\ref{BasicPropertiesOfConvexMinorant}. Let $X=( X_t,t \in[0,1]) $ be a~L\'{e}vy process and consider its convex minorant $C$ on $[0,1]$ as
well as the lower semicontinuous regularization $X^l$ of $X$ given by
$X^l=X\wedge X_-$ [with the convention $X_{-} ( 0) =X_0=0$].

\subsection{\texorpdfstring{Property \protect\ref{ComplementIsLightProperty} of Proposition \protect\ref{BasicPropertiesOfConvexMinorant}}
{Property 1 of Proposition 1}}\label{LightSetProperty}

We will first be concerned with the measure of
\[
\mathscr{P}=\{ t\in[0,1]\dvtx  X^l=C\} .
\]
A first observation is that $\mathscr{P}$ does not vary under changes
in the drift of~$X$.
We now prove that $\mathscr{P}$ has Lebesgue measure zero almost surely.
Indeed, it suffices to see that for each $t\in(0,1)$, $t\notin
\mathscr{P}$
almost surely. If $X$ has unbounded variation, \citet{MR0242261} proves
that
\[
\liminf_{h\to0+}\frac{X_h}{h}=-\infty\qquad\mbox{almost surely (a.s.)}
\]
[see, however, the more recent proof at \citet{MR1875147}], %
and so by the Markov property at each fixed time $t$, we get
\[
\liminf_{h\to0+}\frac{X_{t+h}-X_t}{h}=-\infty\quad\mbox{and}\quad
\limsup_{h\to0+}\frac{X_{t+h}-X_t}{h}=\infty\qquad\mbox{ a.s.}
\]
However, at any $\tau\in\mathscr{P}$, we have
\[
\liminf_{h\to0+}\frac{X_{\tau+h}-X_{\tau}}{h}\geq D ( \tau)
>-\infty\qquad
\mbox{a.s.,}
\]
where $D$ is the right-hand derivative of $C$. If $X$ has bounded
variation, the proof is similar, except that, according to
\citeauthor{MR1406564} [(\citeyear{MR1406564}), Proposition 4, page 81], we get
\[
\lim_{h\to0+}\frac{X_{t+h}-X_t}{h}=d
\]
almost surely, where $d$ is the drift coefficient. We then see that if
$t\in\mathscr{P}\cap(0,1)$, then $D ( t) =d$; the inequality $D ( t)
\leq d$ follows from the preceding display, and by time reversal
we also obtain $d\leq C'_- ( t) $. Taking away the drift, we see that
$t$ then should be a~place where the minimum is achieved. However, $t$
is almost surely not a~time when the minimum is reached:
defining
\[
\tilde X_s=
\cases{
X_{t+s}-X_t,&\quad$\mbox{if } s\leq1-t,$\vspace*{2pt}\cr
X_{1}-X_t+X_{s-(1-t)},&\quad$\mbox{if }1-t\leq s\leq1,$
}
\]
we know that $\tilde X$ has the same law as $X$. Note that the minimum
of $X$ is reached at $t$ if and only if $\tilde X$ remains above zero,
which happens with positive probability only when $0$ is irregular for
$(-\infty,0)$. Hence, $t$ does not belong to $\mathscr{P}$ almost surely
whenever $0$ is regular for $(-\infty,0)$. If this is not the case,
then $0$ is regular for $(0,\infty)$ since $X$ is nonatomic, and
applying same argument to the time reversed process $( X_{( 1-t)
-}-X_1,t\leq1) $ we see then that $t\notin\mathscr{P}$ almost surely
in this remaining case.

\subsection{\texorpdfstring{Property \protect\hyperref[OneJumpPerExcursionProperty]{2} of Proposition \protect\ref{BasicPropertiesOfConvexMinorant}}
{Property 2 of Proposition 1}}\label{OneJumpPerExcursionSubsection}

We will now show that for an excursion interval $( g,d) $ of $X$
above $C$, the jumps of $X$ at $g$ and $d$, denoted $\Delta X_g$
and~$\Delta X_d$, satisfy $\Delta X_g\Delta X_d\geq0$. We first prove
that, thanks to (CD), $X$ does not have jumps of both signs on
the two endpoints of an excursion. The proof depends on different
arguments for bounded and unbounded variation: with unbounded
variation, actually no jumps occur at the endpoints.

If $X$ has unbounded variation, we again use Rogozin's result:
%
\[
\liminf_{h\to0+}\frac{X_h}{h}=-\infty\quad\mbox{and}\quad \limsup
_{h\to0+}\frac{X_h}{h}=\infty,
\]
and adapt Millar's proof of his Proposition 2.4 [\citet
{MR0433606}] %
to see that~$X$ is continuous on $\{ X^l=C\} $. Indeed, for every
$\varepsilon
>0$, let $J_1^\varepsilon,J_2^\varepsilon,\ldots$ be the jumps of
$X$ with size
greater than $\varepsilon$ in absolute value. Then the strong Markov property
applied at $J_i^\varepsilon$ implies that
\[
\liminf_{h\to0+}\frac{X_{J_i^\varepsilon+h}-X_{J_i^\varepsilon
}}{h}=-\infty\quad\mbox{and}\quad \limsup_{h\to0+}\frac{X_{J_i^\varepsilon
+h}-X_{J_i^\varepsilon}}{h}=\infty.
\]
Hence, at any random time $T$ which is almost surely a~jump time of
$X$, we get
\[
\liminf_{h\to0+}\frac{X_{T+h}-X_{T}}{h}=-\infty;
\]
however, if $t\in\{ X^l=C\} $, we see that
\[
\liminf_{h\to0+}\frac{X_{t+h}-X_t}{h}\geq D ( t) >-\infty.
\]

Suppose now that $X$ has bounded variation but infinite L\'{e}vy measure.
Since our problem (jumping to or from the convex minorant) is invariant
under addition of drift, we can assume that the drift coefficient of\vadjust{\goodbreak}
$X$ is zero, and so
\[
\lim_{h\to0+}\frac{X_h}{h}=0
\]
by \citeauthor{MR1406564} (\citeyear{MR1406564}), Proposition 4, page 81. 
We will now prove that almost surely, for every component $( g,d) $
of $\{ C<X^l\} $, we have
%
%
\begin{equation}
\label{JumpsHaveSameSigns}
\Delta X_g\Delta X_d\geq0.
\end{equation}
The argument is similar to the unbounded variation case: at any random
time $T$ which is almost surely a~jump time of $X$, we have
\[
\lim_{h\to0+}\frac{X_{T+h}-X_T}{h}=0.
\]
We deduce that if the slope of $C$ on the interval $( g,d) $ is
strictly positive, then $\Delta X_g\geq0$, and so $X_{g-}=C_g$. By
time-reversal, we see that if the slope of~$C$ is strictly negative on
$( g,d) $, then $\Delta X_d\leq0$ and so $X_d=C_d$. Note that $C$
only has nonzero slopes. Indeed, a~zero slope would mean that the
infimum of $X$ is attained at least twice, a~possibility, that is, ruled
out by Proposition~2.2 of \citet{MR0433606} under assumption
(CD).

\subsection{\texorpdfstring{Property 3 of Proposition \protect\ref{BasicPropertiesOfConvexMinorant}}
{Property 3 of Proposition 1}}

We now see that, almost surely, all excursion intervals of $X$ above
its convex minorant have different slopes. A~different argument is
given for bounded and unbounded variation processes.

When $X$ has unbounded variation on compact sets, let $C^t$ denote the
convex minorant of $X$ on $[0,t]$ so that $C=C^1$. Note that $C^t$ and
$C$ agree up to some random time, which we call $\tau_t$; for every
fixed $t\in(0,1)$, $\tau_t$ cannot equal $t$ as $C_t<X^l_t$ almost
surely, as proved in Section~\ref{LightSetProperty}. We will first
prove that, almost surely, for every $t\in(0,1)\cap\mathbb{Q}$,
whenever the
post~$t$ process touches a~line that extends $C^t$ linearly outwards
from one of the excursion intervals of $C^t$, it crosses it downwards.
To see that this is enough, suppose that there were two excursion
intervals, $(g_1,d_1)$ and $(g_2,d_2)$, with the same associated slope.
Then there would exist $t \in(g_2,d_2)\cap\mathbb{Q}$ such that $g_1
<d_1\leq\tau_t < t$. If the post $t$ process touches the linear
extension of the convex minorant over the interval $(g_1,d_1)$ it must
cross it downwards. This should occur at $d_2$, which contradicts
$C_{d_2} = X^l_{d_2}$.\vspace*{1pt}

To prove the claim that the post $t$ process crosses the extended lines
downwards for each fixed $t\in(0,1)$, let $L_i ( s) =\alpha_i+\beta
_i s$ be the lines extending the segments of $C^t$ (using any ordering
which makes the $\alpha_i$ and $\beta_i$ random variables). Let
\[
T_i=\inf\{ s\geq0\dvtx  X^l_{t+s}-X_t\leq\alpha_i-X_t+\beta_i (t+s)\} .
\]
Hence $T_i$ is a~stopping time for the filtration \mbox{$\mathscr
{F}_{t+s}=\sigma( X_{r}\dvtx r\leq t+s) ,s\geq0$}, with respect to which
$X_{t+s}-X_t,s\geq0$,
is a~L\'{e}vy process. If $X$ jumps below~$L_i$ at time $T_i$, then the
excursion interval of $C$ containing $t$ cannot have\vadjust{\goodbreak} slope~$\beta_i$
(and incidentally,~$\beta_i$ is not a~slope of $C$). Since $X$ has
infinite variation, Rogozin's result quoted above gives
\[
\liminf_{h\to0+}\frac{X_{T_i+h}-X_{t_i}}{h}=-\infty.
\]
Hence, if $X$ is continuous at $T_i$ then $X$ goes below $L_i$
immediately after $T_i$ and $\beta_i$ cannot be a~slope of $C$. We have
seen, however, that in the unbounded variation case, $X$ does not jump
at the vertices of excursion intervals.

When $X$ has bounded variation, the argument is similar except in a~few
places. Suppose the drift of $X$ is zero. We first use
\[
\lim_{h\downarrow0+}\frac{X_{t+h}-X_t}{h}=0
\]
to prove that for every $t\in(0,1)$, whenever the post $t$ process
touches a~linear extension $L_i$ of $C^t$ on an excursion interval
\textit{with positive slope}, by a~jump, it crosses it downwards: this is
clear if $X$ is continuous at $T_i$ or if it jumps into~$L_i$ at $T_i$.
However, $X$ cannot reach $L_i$ from the left and jump away at~$T_i$ by
quasi-continuity of L\'{e}vy processes. By time reversal, we handle the
case of negative slopes, and therefore there are no two excursions
above the convex minorant with the same slope almost surely by the same
arguments as in the unbounded variation case. Again, note that slopes
of $C$ are nonzero since under (CD) the minimum of $X$ is
attained only once by Proposition~2.2 of \citet{MR0433606}.\vspace*{-2pt}

\section{Examples}\label{ExamplesSection}\vspace*{-2pt}

\subsection{L\'{e}vy processes of bounded variation}

Consider a~L\'{e}vy process $X$ with paths of bounded variation on
compact sets and zero drift such that~0 is regular for $(0,\infty)$ but
irregular for $(-\infty,0)$.
Then the cumulative minimum of $X$ is piecewise constant and decreases
by jumps; that is, $X$ reaches a~new minimum by jumping downwards. It
follows that the convex minorant of $X$ on any finite interval has a~finite number of segments of negative slopes until it reaches the
minimum of $X$, and all the excursions above the convex minorant end by
a jump (and begin continuously). However, since the minimum is attained
at a~jump time, say at $\rho$, then $\lim_{t\to0}( X_{\rho
+t}-X_\rho) /t=0$, and since $X_{\rho+\cdot}-X_\rho$ visits
$(0,\infty)$
on any neighborhood of $0$, there cannot be a~segment of the convex
minorant with slope zero, nor a~first segment with positive slope.
Hence $0$ is an accumulation point for positive slopes.

\subsection{The convex minorant of a~Cauchy process}
\mbox{}
\begin{pf*}{Proof of Corollary \protect\ref{CauchyCorollary}}
Let~$X$ be a~symmetric Cauchy process, such that
\[
F ( x) :=\mathbb{P} ( X_1\leq x )=1/2+\arctan ( x) /\pi.
\]
Since
\[
\mathbb{E} ( e^{iu X_t} )= e^{-t|u|},\vspace*{-2pt}\vadjust{\goodbreak}
\]
we see that $X$ is $1$-selfsimilar, which means that $X_t$ has the same
law as $tX_1$ for every $t\geq0$.

If $\Xi_1$ is the point process of lengths and increments of excursions
intervals for the convex minorant on $[0,1]$, its intensity measure
$\nu
_1$ has the following form:
\[
\nu_1 ( dl,dx) =\frac{1}{l} \mathbb{P} ( X_l\in dx ) \,dl.\vspace*{-2pt}
\]
Therefore, the intensity $\tilde\nu_1$ of the point process of lengths
and slopes of excursions intervals for the convex minorant on $[0,1]$,
say $\tilde\Xi_1$, factorizes as
\[
\tilde\nu_1 ( dl,ds) =\frac{1}{l} \mathbb{P} ( X_1\in ds ) \,dl.\vspace*{-2pt}
\]
Let $Y_1,Y_2,\ldots$ be an i.i.d. sequence of Cauchy random variables
independent of~$L$; recall that $F$ is their distribution function.
From the analysis of the point process~$\Xi_1$ in the forthcoming proof
of Lemma~\ref{PPPFromStickBreaking}, the above factorization of the
intensity measure $\tilde\nu_1$ implies that $\tilde\Xi_1$ has the law
of the point process with atoms
%
%
\begin{equation}
\label{AnotherDescriptionForPointProcessInCauchyCase}
\{ ( L_i,Y_i) \dvtx i\geq1\} ;\vspace*{-2pt}
\end{equation}
otherwise said: lengths and slopes are independent for the Cauchy process.

In the converse direction, we see that if lengths and slopes are
independent then~$X$ is a~$1$-selfsimilar L\'{e}vy process. Indeed,
using Theorem~\ref{PointProcessDescriptionForDeterministicAndFiniteHorizon}, we see that
$X_{L_1}/L_1$ and $L_1$ are independent. Let $G$ be the law of
$X_{L_1}/L_1$. Independence of $L_1$ and $X_{L_1}/L_1$ implies that
$X_t/t$ has law $G$ for almost all $t\in(0,1)$, so that $G=F$. As the
law of $X_t/t$ is weakly continuous, we see that $X_t/t$ has law $F$
for all $t\in(0,1)$ and the independence and homogeneity of increments
of $X$ implies that $X_t/t$ has law $F$ for all $t$. However, it is
known that a~$1$-selfsimilar L\'{e}vy process is a~symmetric Cauchy
process, although perhaps seen at a~different speed. See Theorem 14.15
and Example 14.17 of \citet{MR1739520}.

We finish the proof by identifying the law of $( I_x,x\in\mathbb{R}) $.
Informally, $I_x$ is the time in which the convex minorant of $X$ on
$[0,1]$ stops using slopes smaller than~$x$. We then see that $I$ has
the same law as
\[
\tilde I=\Biggl( \sum_{i=1}^\infty L_i{\mathbf{1}}_{Y_i\leq x}, x\in
\mathbb{R}\Biggr) .\vspace*{-2pt}
\]
In contrast, if $U_i,i\geq1$, is an i.i.d. sequence of uniform random
variables on $(0,1)$ independent of $L$, the process $( T_t/T_1,t\in
[0,1]) $ has the representation
\[
\Biggl( \sum_{i=1}^\infty L_i{\mathbf{1}}_{U_i\leq t},t\in[0,1]\Biggr) .\vspace*{-2pt}
\]
With the explicit choice $U_i=F ( Y_i) $, we obtain the result.
\end{pf*}
As a~consequence of Corollary~\ref{CauchyCorollary}, we see that the
set $\mathscr{C}=\{ t\in[0,1]\dvtx  C_t=X_t\wedge X_{t-}\} $ is perfect.\vspace*{-2pt}\vadjust{\goodbreak}

\subsection{The convex minorant of stable processes}
Let $C$ be the convex minorant of the L\'{e}vy process $X$ on $[0,1]$.
We now point out a~dichotomy concerning the set of slopes
\[
\mathscr{S}=\biggl\{ \frac{C_d-C_g}{d-g}\dvtx  ( g,d) \mbox{ is an excursion
interval}\biggr\} ,
\]
when $X$ is a~stable L\'{e}vy process of index $\alpha\in(0,2]$
characterized either by the scaling property
\[
X_{st}\stackrel{d}{=}s^{1/\alpha} X_t,\qquad s>0,
\]
or the following property of its characteristic function:
\[
|\mathbb{E} ( e^{iu X_t} )|=e^{-tc|u|^{\alpha}}.\vspace*{-2pt}
\]

\begin{corollary}
When $\alpha\in(1,2]$, $\mathscr{S}$ has no accumulation points, and
$\mathscr{S}\cap(a,\infty)$ and $\mathscr{S}\cap(-\infty,-a)$ are
almost surely
infinite for all $a>0$. If $\alpha\in(0,1]$, then $\mathscr{S}$ is
dense in
$\mathbb{R}_+$, $\mathbb{R}_-$, or $\mathbb{R}$ depending on if $X$
is a~subordinator, $-X$
is a~subordinator or neither condition holds.\vspace*{-2pt}
\end{corollary}
\begin{pf}
When $\alpha\in(1,2]$, Fourier inversion implies that $X_1$ admits a~continuous and bounded density which is strictly positive. We now make
an intensity measure computation for $a<b$:
\[
\mathbb{E} \bigl( \# \mathscr{S}\cap(a,b) \bigr)
=\int_0^1\int_a^b\frac{1}{t} \mathbb{P} \bigl( X_t/t\in(a,b) \bigr) \,dt.
\]
Using the scaling properties of $X$, we see that near $t=0$, the
integrand is asymptotic to $ct^{-1/\alpha}$ where $c$ is the density of
$X_1$ at zero. Since
\[
\mathbb{E} \bigl( \# \mathscr{S}\cap(a,b) \bigr)<\infty
\]
for all $a<b$, then $\mathscr{S}$ does not contain accumulation points in
$\mathbb{R}$.

If $a>0$, a~similar argument implies that
\[
\mathbb{E} \bigl( \# \mathscr{S}\cap(a,\infty) \bigr)=\infty
\]
since $\mathbb{P} ( X_1>0 )>0$. Unfortunately, this does not imply
that $\#\mathscr{S}\cap(a,\infty)=\infty$ almost surely. However,
from Theorem~\ref{PointProcessDescriptionForDeterministicAndFiniteHorizon}, we see that
$\mathscr{S}\cap[a,\infty)$ has the same law as
\[
\sum_{i\geq1}{\mathbf{1}}_{Y_i\geq aL_i^{1-1/\alpha}},
\]
where $L$ and $Y$ are independent, and $Y_i$ has the same law as $X_1$.
Since $1-1/\alpha>0$ and $L_i\to0$, we see that $Y_i\geq
aL_i^{1-1/\alpha}$ infinitely often, implying that $\#\mathscr{S}\cap
(a,\infty)=\infty$ almost surely.

We have already dealt with the Cauchy case, which corresponds to
$\alpha
=1$, so consider $\alpha\in(0,1)$. Arguing as before, we see that
\[
\# \mathscr{S}\cap(a,b)
\stackrel{d}{=}\sum_{i\geq1}{\mathbf{1}}_{Y_i\in L_i^{1-1/\alpha}(a,b)}.\vadjust{\goodbreak}
\]
Since $1-1/\alpha<0$, we see that $Y_i \in L_i^{1-1/\alpha}(a,b)$
infinitely often as long as $\mathbb{P} ( X_1\in(a,b)>0 )$. Finally, recall
that the support of the law of $X_1$ is $\mathbb{R}_+$, $\mathbb
{R}_-$ or $\mathbb{R}$
depending on if $X$ is a~subordinator, $-X$ is a~subordinator or
neither condition holds.
\end{pf}

\section{Splitting at the minimum and the convex minorant up to an
independent exponential time}
\label{SplittingAtTheMinimumSection}

In this section, we analyze the relationship between Theorem~\ref{PointProcessDescriptionForDeterministicAndFiniteHorizon} and Corollary
\ref{PointProcessDescriptionForIndependentExponentialHorizon} and how
they link with well-known results of the fluctuation theory of L\'{e}vy
processes. We also give a~proof of Corollary~\ref{PointProcessDescriptionForInfiniteHorizon}.

We will first give a~proof of Corollary~\ref{PointProcessDescriptionForIndependentExponentialHorizon} and show how
it leads to a~proof of Theorem~\ref{PointProcessDescriptionForDeterministicAndFiniteHorizon}. While the
implication is based on very well-known results of fluctuation theory,
it is insufficient to prove the more general Theorem~\ref{PointProcessOfExcursions}. Our proof of Theorem~\ref{PointProcessOfExcursions}
is independent of the results of this section.

Let $X$ denote a~L\'{e}vy process with continuous distributions, $C$
its convex minorant on an interval $[0,T]$ (which can be random), $X^l$
the lower-semicontinuous regularization of $X$ given by
$X^l_t=X_t\wedge X_{t-}$ and $\mathscr{O}=\{ s\leq T\dvtx  C_s<X^l_s\} $ is the
open set of excursions from the convex minorant on $[0,T]$. Thanks to
Proposition~\ref{BasicPropertiesOfConvexMinorant} on the basic
properties of the convex minorant, proved in Section~\ref{BasicPropertiesSection}, we see that the point process of lengths and
increments of excursion intervals are equivalently obtained by the
following construction, taken from \citeauthor{MR1739699} [(\citeyear{MR1739699}), Chapter XI]:
define
\[
X^a_t=X_t-at\quad \mbox{and}\quad \underline X^a_t=\min_{s\leq t}X^a_s
\]
as well as
\[
\rho^a=\sup\{ s\leq T\dvtx  X^a_t\wedge X^a_{t-}=\underline X^a_t\}
\quad\mbox{and}\quad
m^a=X^l_{\rho_a}.
\]
The idea behind such definitions is that if $a\mapsto\rho^a$ jumps at
$a$, it is because the convex minorant on $[0,t]$ begins using the
slope $a$ at $\rho^{a-}$ and ends using it at $\rho^a$, while the value
of the convex minorant at the beginning of this interval is $m^{a-}$,
and at the end it is $m^a$. For every fixed $a$, we know that $X^a$
reaches its minimum only once almost surely. However, at a~random $a$
at which $\rho^a$ jumps, the minimum is reached twice, since we know
that slopes are used only once on each excursion interval. From this
analysis, we see that
\[
C_{\rho^a}=X^l_{\rho^a}=m^a
\]
and obtain the following important relationship:
\[
\Xi_T\mbox{ is the point process }\{ ( \rho^a-\rho
^{a-},m^a-m^{a-}) \dvtx  \rho^{a-}<\rho^a\} .
\]
We characterize the two-dimensional process $( \rho,m) $ with the
help of the following results. First of all, according to Millar's
analysis of the behavior of a~L\'{e}vy process at its infimum [cf.
\citet{MR0433606}, Proposition 2.4], 
if $0$ is irregular for $(-\infty, 0)$ then, since $0$ is regular for
$(0,\infty)$, $X^a_{\rho^a}=\underline X^a_{\rho^a}$ almost surely for
each fixed $a$; cf. also the final part of Section~\ref{OneJumpPerExcursionSubsection}.
With this preliminary, Theorem 5 and Lemma~6 from \citeauthor{MR1406564} [(\citeyear{MR1406564}), Chapter
VI] can be written as follows:
\begin{theorem}
\label{FluctuationTheoryTheorem}
Let $T$ be exponential with parameter $\theta$ and independent of~$X$.
For each fixed $a\in\mathbb{R}$, there is independence between the processes
\[
\bigl( X^a_{( t+\rho^a) \wedge T}-m^a,t\geq0\bigr)\quad \mbox{and}\quad
( X^a_{t\wedge\rho^a},t\geq0) .
\]
%
Furthermore,
%
%
\begin{eqnarray}
\label{PechereskiiRogozinFormula}
&&\mathbb{E} \bigl( \exp \bigl( -\alpha\rho^a+\beta( m^a-a \rho^a) \bigr) \bigr)
\nonumber
\\[-8pt]
\\[-8pt]
\nonumber
&&\qquad =\exp
\biggl( {-\int_0^\infty\int_{-\infty}^0 ( 1-e^{-\alpha t+\beta x}) \frac
{e^{-\theta t}}{t} \mathbb{P} ( X_t-at\in dx ) \,dt }\biggr) .
\end{eqnarray}
\end{theorem}

Formula (\ref{PechereskiiRogozinFormula}) was proved initially by
\citet{MR0260005}. Later, \citet{MR588409} showed how
to deduce it by
splitting at the minimum of the trajectory of a~L\'{e}vy process up to
an independent exponential time, a~theme which was retaken by
\citet{MR1406564} to produce the independence assertion of the
previous theorem.
\begin{pf*}{Proof of Corollary
\protect\ref{PointProcessDescriptionForIndependentExponentialHorizon}}
The proof follows [\citet{MR1739699}]. We first show that $(
\rho ,m) $ is a~process with independent increments. Let $a<b$. Note that
$\rho^b-\rho^a$ is the last time that $t$ such that $X_{\rho
^a+t}-m^a-bt$ reaches its minimum, so that Theorem~\ref{FluctuationTheoryTheorem} implies the independence of $\rho
^{a+b}-\rho
^a$ and $\sigma( X_{\cdot\wedge\rho^a}) $; denote the latter $\sigma$-field as $\mathscr{F}
^a$. Also, note that $m^b-m^a$ is the minimum of $X_{( \rho ^a+t)
\wedge T}-m^a-bt, t\geq0$. Hence there is also independence
between $m^b-m^a$ and $\mathscr{F}^a$. Finally, note that if $a'\leq
a,$ $( \rho^{a'},m^{a'}) $ are $\mathscr{F}^a$ measurable since $\rho
^{a'}$ is the last
time that $X_{\cdot\wedge\rho^a}-a'\cdot$ reaches its minimum on
$[0,\rho^a]$, and $m^{a'}$ is the value of this minimum.

From the above paragraph, we see that the point process of jumps of
$( \rho,m) $, that is, $\Xi$, is a~Poisson random measure: this
would follow from (a~bidimensional extension of) Theorem 2 and
Corollary 2 in Gihman and Skorohod [(\citeyear{MR0375463}), Chapter~IV.1, pages 263--266] which
affirm that
the jump process of a~stochastically continuous process with
independent increments on $\mathbb{R}_+$ is a~Poisson random measure
on $\mathbb{R}
_+\times\mathbb{R}_+$. To show that $( \rho,m) $ is stochastically
continuous, we show that it has no fixed discontinuities; this follows
because for every fixed $a\in\mathbb{R}$, the minimum of $X^a$ is
reached at
an unique point almost surely, which implies that, for every fixed $a$,
almost surely, neither $\rho$ nor $m$ can jump at $a$. To compute the
intensity measure $\nu$ of $\Xi_{T}$, note that the pair $( \rho
^a,m^a) $ can be obtained from $\Xi_T$ as
%
%
\begin{equation}
\label{RhoMInTermsOfXiEquation}
( \rho^a,m^a) =\mathop{\sum_{( u,v) \in\mathscr{I}}}_{ C_v-C_u\leq
a( v-u) }( v-u,C_v-C_u) .
\end{equation}
The above equality contains the nontrivial assertion that the additive
process $( \rho,m) $ has no deterministic component or, stated
differently, that it is the sum of its jumps. For the process $\rho$,
this follows because
\[
\mathop{\sum_{( u,v) \in\mathscr{I}}}_{C_v-C_u\leq a( v-u) }(
v-u) =\mbox{Leb} ( \mathscr{O}\cap\{ t\leq T\dvtx  C'_t\leq a\} )
\]
which, since $\mbox{Leb} ( \mathscr{O}) =T$ and $C'$ is
nondecreasing, gives
\[
\mathop{\sum_{( u,v) \in\mathscr{I}}}_{ C_v-C_u\leq a( v-u) }(
v-u) =\sup\{ t\leq T\dvtx C'_t\leq a\} =\rho^a.
\]
To discuss the absence of drift from $m$, let $m^C$ be the signed
measure which assigns each interval $( u,v) $ the quantity
$C_v-C_u$. (Because $C'$ is nondecreasing, it is trivial to prove the
existence of such a~signed measure, to give a~Hahn decomposition of it
and to see that it is absolutely continuous with respect to Lebesgue
measure.) Then
\begin{eqnarray*}
\mathop{\sum_{( u,v) \in\mathscr{I}}}_{ C_v-C_u\leq a( v-u) }( C_v-C_u)
&=&m^C ( \mathscr{O}\cap\{ t\leq T\dvtx  C'_t\leq a\} )
\\
&=&m^C ( \{ t\leq T\dvtx  C'_t\leq a\} )
=C_{\rho^a}=X^l_{\rho^a}
=m^a.
\end{eqnarray*}

From (\ref{RhoMInTermsOfXiEquation}), we get
\[
\mathbb{E} \bigl( \exp ( -\alpha\rho^a+\beta m^a) \bigr)=\exp \biggl( -\int
_0^\infty\int _{-\infty}^{at}( 1-e^{-\alpha t+\beta x}) \nu ( dt,
dx) \biggr) ,
\]
while from the Pe{\v{c}}erski{\u\i}--Rogozin formula (\ref{PechereskiiRogozinFormula}), we obtain
\begin{eqnarray*}
&&\mathbb{E} \bigl( \exp ( -\alpha\rho^a+\beta m^a) \bigr)
\\
&&\qquad=\exp \biggl( -\int_0^{\infty}\int_{-\infty}^0 \bigl( 1-e^{-( \alpha
-a\beta) t+\beta x}\bigr) \frac{e^{-\theta t}}{t} \mathbb{P} ( X_t-at\in
dx ) \,dt\biggr)
\\
&&\qquad=\exp \biggl( -\int_0^{\infty}\int_{-\infty}^{at} ( 1-e^{-\alpha
t+\beta x}) \frac{e^{-\theta t}}{t} \mathbb{P} ( X_t\in dx ) \,dt\biggr)
\end{eqnarray*}
giving
\[
\nu ( dt,dx) =\frac{e^{-\theta t}}{t} \mathbb{P} ( X_t\in dx ) \,dt.
\]
\upqed\end{pf*}

We now remark on the equivalence between Theorem~\ref{PointProcessDescriptionForDeterministicAndFiniteHorizon} and Corollary
\ref{PointProcessDescriptionForIndependentExponentialHorizon} and how
either of them implies Corollary~\ref{PointProcessDescriptionForInfiniteHorizon}.

Let $L$ be an uniform stick-breaking sequence and $X$ a~L\'{e}vy
process with continuous distributions which are independent. Let $S$
be\vadjust{\goodbreak}
the partial sum sequence associated to $L$, and consider the point
process $\tilde\Xi_t$ with atoms at
\[
\{ ( tL_i, X_{t S_i}-X_{tS_{i-1}}) \} .\vspace*{-2pt}
\]

\begin{lemma}
\label{PPPFromStickBreaking}
If $T$ an exponential random variable of parameter $\theta$ independent
of $( X,L) $, $\tilde\Xi_T$ is a~Poisson point process with intensity
%
%
\begin{equation}
\label{IntensityInPPPFromStickBreaking}
\mu_\theta ( dt,dx) =\frac{e^{-\theta t}}{t} \,dt\,\mathbb{P} ( X_t\in
dx ).
\end{equation}
\end{lemma}
\begin{pf}
We recall the relationship between the Gamma subordinator and the
stick-breaking process, which was found by McCloskey in his unpublished
PhD thesis [\citet{McCloskey}] and further examined and extended by
\citet{MR1156448}. Recall that a~Gamma process is a~subordinator
$( \Gamma_t,t\geq0) $ characterized by the Laplace exponent
\[
\mathbb{E} ( e^{-q \Gamma_t} )=\biggl( \frac{\theta}{\theta+q}\biggr) ^t=\exp
\biggl( -t\int_0^\infty( 1-e^{-q x}) \frac{e^{-\theta x}}{x} \,dx\biggr) ;
\]
the law of $\Gamma_1$ is exponential of parameter $\theta$.
It is well known that $( \Gamma_t/\Gamma_1, t\leq1) $ is
independent of $\Gamma_1$. Also, it was proved [\citet{McCloskey}; \citet{MR1156448}] that the size-biased permutation of the jumps
of $( \Gamma_t/\Gamma_1, t\in[0,1]) $ has the same law as the stick-breaking
process on $[0,1]$. Hence if $L$ is a~stick-breaking process
independent of the exponential $T$ of parameter $\theta$, then the
point process with atoms at $\{ TL_1,TL_2,\ldots\} $ has the same law
as the point process with atoms at the jumps of a~Gamma subordinator
(of parameter $\theta$) on $[0,1]$ or, equivalently, a~Poisson point
process with intensity $e^{-\theta x}/x \,dx$.

If $S$ is the partial sum sequence associated to $L$, conditionally on
\mbox{$T=t$} and $L=( l_1,l_2,\ldots) $, $( X_{TS_i}-X_{TS_{i-1}},i\leq1) $
are independent and the law of
$X_{TS_i}-X_{TS_{i-1}}$ is that of $X_{t l_i}$. We deduce that the
point process with atoms $\{ ( TL_i,X_{TS_i}-X_{TS_{i-1}}) ,i\geq 1\}
$ is a~Poisson point process with the intensity~$\mu_\theta$ of (\ref{IntensityInPPPFromStickBreaking}), as shown, for example, in
\citeauthor{MR1876169} [(\citeyear{MR1876169}), Proposition~12.3, page~228] using the notion of
randomization of point processes.\vspace*{-2pt}
\end{pf}

Lemma~\ref{PPPFromStickBreaking} shows how Theorem~\ref{PointProcessDescriptionForDeterministicAndFiniteHorizon} implies
Corollary~\ref{PointProcessDescriptionForIndependentExponentialHorizon}.

Conversely, if we assume Corollary~\ref{PointProcessDescriptionForIndependentExponentialHorizon}, we know that
$\tilde\Xi_T$ has the same law as the point process of lengths and
increments of excursions intervals on the interval $[0,T]$. However, if
$ \Xi_t$ is the point process of lengths and increments of excursion
intervals on $[0,t]$, then
\[
\int_0^\infty\theta e^{-\theta t}\mathbb{E} ( e^{-\Xi_t f} )
\,dt=\mathbb{E} ( e^{-\tilde\Xi_T f} )=\int_0^\infty\theta
e^{-\theta t} \mathbb{E} ( e^{-\tilde\Xi_t f} ) \,dt
\]
which implies that
\[
\mathbb{E} ( e^{- \Xi_t f} )=\mathbb{E} ( e^{-\tilde\Xi_t f} )
\]
for continuous and nonnegative $f$. However, this implies the identity
in law between $\Xi_t$ and $\tilde\Xi_t$, giving
Theorem~\ref{PointProcessDescriptionForDeterministicAndFiniteHorizon}.\vadjust{\goodbreak}

Let us pass to the proof of Corollary~\ref{PointProcessDescriptionForInfiniteHorizon}. Abramson and Pitman show
the discrete time analog using a~Poisson thinning procedure.\vspace*{-2pt}
\begin{pf*}{Proof of Corollary
\protect\ref{PointProcessDescriptionForInfiniteHorizon}}
Suppose $l=\liminf_{t\to\infty}X_t/t\in(-\infty,\infty]$. Then there
exists $a\in\mathbb{R}$ and $T>0$ such that $X_t>at$ for all $t>T$.
If $C^T$
is the convex minorant of $X$ on $[0,T]$, and $\rho^a$ is the first
instant at which the derivative of $C^T$ is greater than $a$, then the
convex function
\[
\tilde C_t=
\cases{
C^T,&\quad$\mbox{if }t< \rho^a,$\vspace*{2pt}\cr
C^T_{T_a}+a( t-T), &\quad$\mbox{if }t\geq\rho^a,$}
\]
lies below the path of $X$ on $[0,\infty)$, implying $C^\infty$, the
convex minorant of $X$ on $[0,\infty)$, is finite for every point of
$[0,\infty)$.

Conversely, if $C^\infty$ is finite on $[0,\infty)$, for any $t>0$ we
can let $a=\break \lim_{h\to0+}( C_{t+h}-C_t) /h\in\mathbb{R}$ and note that
$\liminf_{s\to\infty}X_s/s\geq a$.

From \citet{MR0336806} we see that, actually, $\lim_{t\to
\infty}X_t/t$
exists and it is finite if and only if $\mathbb{E} ( |X_1| )<\infty$ and
$\mathbb{E} ( X_1 )=l$.
Note that the right-hand derivative of~$C^\infty$ is never strictly
greater than $l$.
This derivative cannot equal $l$: if $l=\infty$ this is clear while if
$l<\infty$, it follows from the fact that the zero mean L\'{e}vy
process $X_t-l t$ visits $(-\infty,0)$ [as can be proved, e.g., by
embedding a~random walk and using, for example, by \citet{MR0040610}; \citet{MR0133148}].
However, the derivative also surpasses any level $a<l$. This follows
from the definitions of $l$ and $C^\infty$: if the derivative of
$C^\infty$ were always less than $l-\varepsilon$, since $X_t$
eventually stays
above every line of slope $l-\varepsilon/2$, we would be able to
construct a~convex function greater than $C^\infty$ and below the path of $X$.

If $a<l$, let $L_a$ be the last time the derivative of $C^\infty$ is
smaller than $a$. Then for $t>L_a$, we see that
\[
C^{L_a}=C^t=C^\infty\qquad\mbox{on $[0,L_a]$}.
\]
We will now work with $C^{T_\theta}$, where $T_\theta$ is exponential
of parameter $\theta$ and independent of $X$. Then on the set $\{
L_a<T_\theta\} $, which has probability tending to $1$ as $\theta\to0$,
we have $C^{L_a}=C^{T_\theta}=C^{\infty}$ on $[0,L_a]$. Recall,
however, that if~$\Xi_\theta$ is a~Poisson point process with intensity
\[
\mu_\theta ( dt,dx) =\frac{e^{-\theta t}}{t} \mathbb{P} ( X_t\in dx
) \,dt,
\]
then $\Xi_\theta$ has the law of the lengths and increments of
excursions of $X$ above~$C^{T_\theta}$ by Corollary~\ref{PointProcessDescriptionForIndependentExponentialHorizon}. We deduce
that for every $a<l$ the restriction of $\Xi_\theta$ to $\{ ( t,x) \dvtx
x<at\} $ converges in law as $\theta\to0$ to the point process
with atoms at the lengths and increments of excursions of $X$ above
$C^\infty$ with slope less than $a$. Hence, the excursions of $X$ above
$C^\infty$ with slopes $<a$ form a~Poisson point process with intensity
\[
\frac{{\mathbf{1}}_{x<at}}{t} \mathbb{P} ( X_t\in dx ) \,dt.
\]
It suffices then to increase $a$ to $l$ to obtain the stated
description of $\Xi_\infty$.\vadjust{\goodbreak}
\end{pf*}
%

Basic to the analysis of this section has been the independence result
for the pre and post minimum processes up to an independent exponential
time as well as the Pe{\v{c}}erski{\u\i} and Rogozin formula stated in
Theorem~\ref{FluctuationTheoryTheorem}. Theorem~\ref{FluctuationTheoryTheorem} is the building block for the fluctuation
theory presented in \citeauthor{MR1406564} [(\citeyear{MR1406564}), Chapter VI] and is
obtained there
using the local time for the L\'{e}vy process reflected at its
cumulative minimum process. In the following sections, we will reobtain
Theorem~\ref{PointProcessDescriptionForDeterministicAndFiniteHorizon}
and Corollaries~\ref{PointProcessDescriptionForIndependentExponentialHorizon} and~\ref{PointProcessDescriptionForInfiniteHorizon} appealing only to the basic
results of the convex minorant of Section~\ref{BasicPropertiesSection}
(and without the use of local time). In particular, this implies the
first part of Theorem~\ref{FluctuationTheoryTheorem}, from which the
full theorem follows as shown by \citet{MR1406564}. 
Indeed, assuming Theorem~\ref{PointProcessOfExcursions}, if $T$ is
exponential with parameter $\theta$ and independent of $X$, if $\rho^a$
is the last time $X^l_t-a t$ reaches its minimum on $[0,T]$, and $m^a$
is the value of this minimum, we see that
\[
\bigl( X^a_{( t+\rho^a) \wedge T}-m^a, t\geq0\bigr)
\]
can be obtained from the Poisson point process of excursions of $X$
above its convex minorant with slopes $>a$, while
\[
( X^a_{t\wedge\rho^a}, t\geq0)
\]
is obtained from the excursions with slopes $\leq a$. Since the process
of excursions (up to an independent time) is a Poisson point process,
we obtain the independence of the pre and post minimum processes.

Here is another example of how the description of the convex minorant
up to an independent exponential time leads to a~basic result in
fluctuation theory: according to Rogozin's criterion for regularity of
half-lines, $0$ is irregular for $(0,\infty)$ if and only if
%
%
\begin{equation}
\label{RogozinsCriterion}
\int_0^1 \mathbb{P} ( X_t\leq0 )/t \,dt<\infty.
\end{equation}
To see how this might be obtained from Corollary~\ref{PointProcessDescriptionForIndependentExponentialHorizon}, we note that
the probability that $X$ does not visit $(0,\infty)$ on some
$(0,\varepsilon)$
is positive if and only if the convex minorant up to $T_\theta$ has
positive probability of not having negative slopes. By Theorem~\ref{PointProcessDescriptionForIndependentExponentialHorizon} this happens
if and only if
\[
\int_0^\infty\mathbb{P} ( X_t\leq0 )e^{-\theta t}/t \,dt<\infty,
\]
which is of course equivalent to (\ref{RogozinsCriterion}).

\section{Absolute continuity of the minimum and its location}
\label{ApplicationsSection}
\mbox{}
\begin{pf*}{Proof of Theorem \protect\ref{AbsoluteContinuityTheorem}}
Since $0$ is regular for both half-lines, the L\'{e}vy process $X$
satisfies assumption (CD), and we can apply Theorem~\ref{PointProcessDescriptionForDeterministicAndFiniteHorizon}.

Let $L$ be an uniform stick-breaking process independent of $X$, and
define its partial sum and residual processes $S$ and $R$ by
\[
S_0=0,\qquad S_{i+1}=S_i+L_{i+1}\quad \mbox{and}\quad R_i=1-S_i.\vadjust{\goodbreak}
\]
Set
\[
\Delta_i=X_{S_i}-X_{S_{i-1}}.
\]
Then the time of the minimum of the L\'{e}vy process $X$ on $[0,1]$,
has the same law as
\[
\rho=\sum_{i=1}^\infty L_i{\mathbf{1}}_{\Delta_i<0},
\]
while the minimum of $X$ on $[0,1]$ (denoted $\underline X_1$) and
$X_1-\underline X_1$ have the same laws as
\[
\sum_{i=1}^\infty\Delta_i{\mathbf{1}}_{\Delta_i<0} \quad\mbox{and}\quad
\sum
_{i=1}^\infty\Delta_i{\mathbf{1}}_{\Delta_i>0}.
\]
The basic idea of the proof is to decompose these sums at a~random
index~$J$; in the case of $\rho$,
into
%
%
\begin{equation}
\label{DecompositionForMinimumAtRandomIndex}
\Sigma_J=\sum_{i=1}^JL_i{\mathbf{1}}_{\Delta_i<0} \quad\mbox{and}\quad
\Sigma
^J=\sum_{i=J+1}^\infty L_i{\mathbf{1}}_{\Delta_i<0}.
\end{equation}
The random index (actually a~stopping time for the sequence $\Delta$)
is chosen so that $\Sigma_J$ and $\Sigma^J$ are both positive, and
$( R_J,\Sigma_J) $ has a~joint density, which is used to provide a~density
for $\Sigma$ using the conditional independence between~$\Sigma_J$ and $\Sigma^J$ given $R_J$.

Let $I$ be any stopping time for the sequence $\Delta$ which is finite
almost surely. We first assert that the sequence $( \Delta _{I+i-1})
_{i\geq1}$ has both nonnegative and strictly negative terms
if 0 is regular for both half-lines. Indeed, if $0$ is regular for
$(-\infty,0)$, this implies that the convex minorant of $X$ has a~segment of negative slope almost surely, which implies the existence of
$i$ such that $\Delta_i<0$ almost surely. If $0$ is regular for
$(0,\infty)$, a~time-reversal assertion proves also the existence of
nonnegative terms in the sequence $\Delta$. On the other hand,
conditionally on $I=i$ and $L_1=l_1,\ldots,L_i=l_i$, the sequence
$( \Delta_{i-1+j},j\geq1) $ has the same law as the sequence
$\Delta$ but obtained from the L\'{e}vy process $X_{( 1-l_1-\cdots
-l_i) t,t\geq0}$ which shares the same regularity as $X$, which implies
the assertion.
\begin{longlist}[1.]
\item[1.]
Let $I$ and $J$ be defined by
\[
I=\min\{ i\geq1\dvtx  \Delta_i\geq0\} \quad\mbox{and}\quad J=\min\{ j\geq
I\dvtx \Delta_j<0\} .
\]
By the preceding paragraph, we see that $I$ and $J$ are both finite
almost surely. Hence, the two sums $\Sigma_J$ and $\Sigma^J$ of (\ref{DecompositionForMinimumAtRandomIndex}) are both in the interval $(0,1)$
and we have
\[
\rho=\Sigma_J+\Sigma^J.
\]

We now let
\[
f ( t) =\mathbb{P} ( X_t\leq0 )\vadjust{\goodbreak}
\]
which will allow us to write the density of $( \Sigma_{J},R_{J}) $;
this follows from the computation
\begin{eqnarray*}
&&\mathbb{P} ( J=j,L_1\in dl_1,\ldots,L_{j}\in dl_{j} )
\\
&&\qquad =\sum_{i=1}^{j-1} \prod_{k<i}f ( l_k) \prod_{i\leq k<j}\bigl( 1-f (
l_k) \bigr) f ( l_j) \mathbb{P} ( L_1\in l_1,\ldots, L_j\in l_j )
\end{eqnarray*}
valid for $j\geq2$. For $2\leq i<j$, let
\[
g_{i,j} ( l_1,\ldots,l_j) =( l_1,\ldots,l_{i-1},l_i,\ldots
,l_{j-2},l_1+\cdots+l_i+l_j,1-l_i-\cdots-l_j) ,
\]
and define
\[
g_{1,2} ( l_1,l_2) =( l_2,1-l_2-l_2)
\]
as well as
\[
g_{1,j} ( l_1,\ldots,l_j) =( l_1,\ldots ,l_{j-2},l_j,1-l_1-\cdots -l_j)
\]
for $j\geq3$. Then $g_{i,j}$ is an invertible linear transformation on
$\mathbb{R}^j$, and so if~$B$ is a~Borel subset of $\mathbb{R}^j$ of
Lebesgue measure
zero, then $g_{i,j}^{-1} ( B) $ also has Lebesgue measure zero. If $A$ is
a Borel subset of $\mathbb{R}^2$ with Lebesgue measure zero, we get
\[
\mathbb{P} \bigl( ( \Sigma_J,R_J) \in A \bigr)\leq\sum_{j=2}^\infty\sum_{i=1}^{j-2}
\mathbb{P}\bigl ( ( L_1,\ldots,L_j) \in g_{i,j}^{-1} ( \mathbb
{R}^{j-2}\times A) \bigr)=0.
\]
Hence, there exists a~function $g$ which serves as a~joint density of
$( \Sigma_J,R_J) $. We can then let
\[
g_r ( l) =\frac{g ( l,r) }{\int g ( l',r) \, dl'}
\]
be a~version of the conditional density of $\Sigma_J$ given $R_J=r$.

Using the construction of the stick breaking process and the
independence of increments of $X$ we deduce that
\[
\tilde L=\biggl( \frac{L_{J+i}}{R_{J}},i\geq1\biggr)
\]
is independent of $( L_{i\wedge J},\Delta_{i\wedge J},i\geq1) $
and has the same law as $L$. Furthermore, the sequence $( \Delta
_{J+i},i\geq1) $ is conditionally independent of $( L_{i\wedge J},
\Delta_{i\wedge J}) $ given~$R_J$.

We therefore obtain the decomposition
\[
\rho=\Sigma_{J}+R_{J}\rho^{J},
\]
where
\[
\rho^{J}=\sum_{i=1}^{\infty} \frac{L_{i+J}}{R_{J}}{\mathbf
{1}}_{\Delta _{i+J}<0}=\frac{\Sigma^J}{R_J}.\vadjust{\goodbreak}
\]
Since $\rho^{J}$ is a~function of $\tilde L$, $( \Delta _{J+i},i\geq
1) $ and $R_J$, then $\rho^{J}$ and $\Sigma_{J}$ are
conditionally independent given $R_{J}$. Hence $g_{R_{J}}$ is also a~version of the conditional density of $\Sigma_{J}$ given $R_{J}$ and
$\rho^{J}$, and we can then write
%
%
\begin{equation}
\label{DensityOfMinimum}
\mathbb{P} ( \rho\in dt )=dt \int g_{r} ( t-ry) \mathbb{P} (
R_{J}\in dr,\rho ^{J}\in dy )
\end{equation}
on $\{ J<\infty\} $.

Finally, it remains to see that the density for $\rho$ displayed in
equation (\ref{DensityOfMinimum}) is positive on $(0,1)$. We remark
that the density of $( R_J,\Sigma_J) $ is positive on
\[
\{ ( r,\sigma) \dvtx 0<\sigma<1-r<1\} .
\]
Indeed, taking $r,\sigma$ as in the preceding display, we have the
explicit computation
\begin{eqnarray*}
&&\mathbb{P} ( J=2,\Sigma_J\in d\sigma, R_J\in dr )\\
&&\qquad=\mathbb{P} ( \Delta_1>0,\Delta_2<0, L_2\in d\sigma,1-L_1-L_2\in
dr )
\\
&&\qquad=\bigl( 1-f ( 1-\sigma-r) \bigr) f ( \sigma) {\mathbf{1}}_{0<\sigma
<1-r<1}\frac{1}{1-\sigma-r} \,dr \,d\sigma.
\end{eqnarray*}
On the other hand, given $t\in(0,1)$, $\mathbb{P} ( R_J<1-t )>0$. Indeed,
\begin{eqnarray*}
\mathbb{P} ( R_J<1-t )
&\geq& \mathbb{P} ( R_J<1-t,J=2 )
\\
&=&\int\int\mathbb{P} ( \Delta_1\geq0,\Delta_2<0, L_1\in
dl_1,1-L_1-L_2\in dl_2 ){\mathbf{1}}_{l_2\leq1-t}
\\
&=&\int\int\bigl( 1-f ( l_1) \bigr) f ( 1-l_1-l_2) \frac
{1}{1-l_1}{\mathbf{1}}_{0<l_2<1-l_1}{\mathbf{1}}_{l_2<1-t}
\\
&>&0,
\end{eqnarray*}
since $f$ and $1-f$ are strictly positive on $(0,1)$ since $0$ is
regular for both half lines and so the support of the law of $X_t$ is
$\mathbb{R}$ for all $t>0$. Going back to equation (\ref{DensityOfMinimum}),
we see that, given $t\in(0,1)$, on the set $\{ ( r,y) \dvtx  0<r<1-t\} $ we
have $t-ry<t<1-r$, and so the density $g_r ( t-ry) $ is
positive. Hence the integral in equation~(\ref{DensityOfMinimum}) is positive.

\item[2.] The proof of absolute continuity of the time and value of the
minimum of $X$ on $[0,1]$ is similar, except that further hypotheses
are needed.

First, the value of the minimum of $X$ on $[0,1]$ has the same
distribution as
\[
m:=\sum_{i=1}^\infty\Delta_i{\mathbf{1}}_{\Delta_i\leq0}.
\]
Since the law of $X_t$ is absolutely continuous with respect to
Lebesgue measure for all $t>0$, we have
\[
( \rho,m) =( \Sigma_J,m_J) +( R_J\rho^J,m^J),
\]
where
\[
m_J=\sum_{i\leq J}\Delta_i{\mathbf{1}}_{\Delta_i<0}
\quad\mbox{and}\quad
m^J=\sum_{i=1}^{\infty}\Delta_{J+i}{\mathbf{1}}_{\Delta_{J+i}>0}.
\]
We now prove that:
\begin{longlist}[(a)]
\item[(a)]$( \rho_J,m_J) $ has a~conditional density with respect to $R_J$;
\item[(b)]$( \rho_J,m_J) $ and $( \rho^J,m^J) $ are conditionally
independent given $R_J$.
\end{longlist}
The second assertion follows from our previous analysis of conditional
independence in the sequences $L$ and $\Delta$. The first assertion
follows from the fact that $( \Sigma_J,R_J,\Delta_J) $ admit a~density on $\{ J=j\} $, by a~computation similar to the one for $(
\Sigma_J,R_J) $
\begin{eqnarray*}
&&\mathbb{P} ( J=j, L_1\in dl_1,\ldots,L_{j}\in dl_{j}, \Delta_1\in
dx_1,\ldots ,\Delta_{j}\in dx_{j} )
\\
&&\qquad=\sum_{i=1}^{j-1}{\mathbf{1}}_{x_1.\ldots,x_{i-1}<0,x_i,\ldots
,x_{j-1}>0,x_j<0}\mathbb{P} ( X_{l_1}\in dx_1 )\cdots\mathbb{P} (
X_{l_{j}}\in dx_{j} )\\
&&\hspace*{26pt}\qquad{}\times \mathbb{P} ( L_1\in dl_1,\ldots,L_{i}\in dl_{j} )
\end{eqnarray*}
so that on $\{ J=j\} $, $( L_1,\ldots,L_{J},\Delta_1,\ldots,\Delta
_{J}) $ admit a~density with respect to Lebesgue measure, and since
$( \Delta_J,R_{J},\Delta_{J}) $ is the image under a~surjective
linear map of the former variables, the latter admit a~joint density.
Let $f_r$ be a~version of the conditional density of $( \Sigma
_J,\Delta_J) $ given $R_J=r$. We then get
%
%
\begin{eqnarray}
\label{JointDensityForRhoM}
&&\mathbb{P} ( \rho\in dt, m\in dx )
\nonumber
\\[-8pt]
\\[-8pt]
\nonumber
&&\qquad=dt \,dx \int f_{r} ( t-rs,x-y)
\mathbb{P} ( \rho^J\in ds, m^J\in dy, R_J\in dr ).
\end{eqnarray}
Regarding the equivalence of the law of $( \rho,m) $ and Lebesgue
measure on $(0,1)\times(-\infty,0)$, note that a~version of the
density of $( R_J, \rho_J,m_J) $ is positive on
\[
\{ ( r,s,x) \dvtx 0\leq r+s\leq1, x<0\} .
\]
Indeed, we have, for example,
\begin{eqnarray*}
&&\mathbb{P} ( \Delta_1<0,\Delta_2>0,R_I\in dr,\Sigma_I\in ds,m_I\in
dx )
\\
&&\qquad =\mathbb{P} ( X_s\in dx )\bigl( 1-f ( 1-r-s) \bigr) \frac{1}{1-s}{\mathbf
{1}}_{0\leq r+s\leq1}{\mathbf{1}}_{x\leq0}.
\end{eqnarray*}
Since the law of $( \rho^J,m^J,R_J)$, by analogy with the case of
$\rho$, is seen to charge the set $\{ ( s,y,r) \dvtx t<1-r, x<y\} $, we
conclude that the expression for the joint density of $( \rho,m) $
given in equation (\ref{JointDensityForRhoM}) is strictly positive.
\item[3.] The proof of the absolute continuity of $( \underline
X_1,X_1-\underline X_1) $ follows the same method of proof, starting
with the fact that these random variables have\vadjust{\goodbreak} the same joint law as
\[
( \Delta^-,\Delta^+) =\sum_{i=1}^\infty\Delta_i( {\mathbf
{1}}_{\Delta _i<0},{\mathbf{1}}_{\Delta_i>0}) ,
\]
which we can again decompose at the random index
\[
I=\min\{ i\geq1\dvtx  \mbox{there exist $j,j\leq i$ such that $\Delta
_j<0$, $\Delta_{j'}>0$}\}
\]
into
\[
( \Delta^-,\Delta^+) =( \Delta^-_I,\Delta^+_I) +( \Delta
^{-,I},\Delta^{+,I}),
\]
where
\[
( \Delta^-_I,\Delta^-_I) =\sum_{i\leq I}\Delta_i( {\mathbf
{1}}_{\Delta _i<0},{\mathbf{1}}_{\Delta_i>0}) .
\]
Since:
\begin{longlist}[(a)]
\item[(a)]$( R_I,\Delta^-_I,\Delta^+_I) $ have a~joint density which
can be taken positive on $(0,1)\times(-\infty,0)\times(0,\infty)$, and
\item[(b)]$( \Delta^-_I,\Delta^+_I) $ and $( \Delta^{-,I},\Delta
^{+,I}) $ are conditionally independent given $R_I$,
\end{longlist}
we see that $( \Delta^-,\Delta^+) $ admit a~joint density which can
be taken positive on $(-\infty,0)\times(0,\infty)$.
\end{longlist}
\upqed\end{pf*}

\section{An invariant path transformation for L\'{e}vy processes}
\label{InvariantPathTransformationSection}

The aim of this section is to prove Theorem~\ref{InvarianceUnderPathTransformation}. This will be done (almost) by
applying the continuous mapping theorem to the embedded random walk
$( X_{k/n},k=0,\ldots, n) $ and a~continuous function on Skorohod
space. The argument's technicalities are better isolated by focusing
first on some special cases in which the main idea stands out.
Therefore, we first comment on the case when $X$ has continuous sample
paths, then we handle the case when~$X$ has paths of unbounded
variation on compact intervals, to finally settle the general case.

We rely on a~discrete version of Theorem~\ref{InvarianceUnderPathTransformation}, which was discovered by
\citet{2010AbramsonPitman}. Let $S^n=( S^n_t,t\in[0,n]) $ be the process
obtained by interpolating between the values of $n$ steps of a~random
walks which jumps every $1/n$, and let $C^n$ be its convex minorant.
Let $V^n_0,V^n_1,\ldots, V^n_k$ be the endpoints of the segments
defining the convex minorant $C^n$. Let $U_n$ be uniform on $\{
1/n,\ldots,1\} $. Since there exists an unique $j$ such that
\[
U_n\in(V^n_j,V^n_{j+1}],
\]
let us define
\[
g_n=V^n_j\quad \mbox{and}\quad d_n=V^n_{j+1}
\]
as the excursion interval of $S^n$ above $C^n$ which straddles $U_n$.
Mimicking the definition of the path transformation (\ref{PathTransformationForLevyProcess}), let us define
\[
S^{n,U_n}_t=
\cases{
S^n_{U_n+t}-S^n_{U_n},&\quad$\mbox{if } 0\leq t\leq
d_n-U_n,$\vspace*{2pt}\cr
S^n_{d_n}-S^n_{U_n}+S^n_{g_n+t-(d_n-U_n)}-S^n_{g_n},&\quad$\mbox{if } d_n-U_n\leq t\leq d_n-g_n,$\vspace*{2pt}\cr
S^n_{d_n}+S^n_{t-( d_n-g_n) },&\quad $\mbox{if } d_n-g_n\leq t\leq
d_n,$\vspace*{2pt}\cr
S^n_{t},&\quad $\mbox{if } d_n\leq t.$}
\]

\begin{theorem}[{[\citet{2010AbramsonPitman}]}]
\label{2010AbramsonPitmanTheorem}
If the distribution function of $S^n_{1/n}$ is continuous, then the pairs
\[
( U_n,S^n)\quad \mbox{and} \quad( d_n-g_n,S^{n,U_n})
\]
have the same law.
\end{theorem}

To prove Theorem~\ref{InvarianceUnderPathTransformation} we will use
Theorem~\ref{2010AbramsonPitmanTheorem} with the random walk obtained
by sampling our L\'{e}vy process $X$ at points of the form $1/n$ and
take the limit as $n\to\infty$. The details are a~bit technical in
general but simplify considerably when~$X$ is continuous or when it
reaches its convex minorant continuously.

The main tool for
the passage to the limit
is a~lemma regarding approximation of the endpoints of the interval of
the convex minorant that contains a~given point.
Let $f\dvtx [0,1]\to\mathbb{R}$ be a~c\`{a}dl\`{a}g\ function which starts
at zero
and is left continuous at $1$ and $c$ its convex minorant.
Let also $f^l=f\wedge f_-$ be the lower semicontinuous regularization
of $f$, 
and define with it the excursion set away from the convex minorant
$\mathscr{O}=\{ c<f^l\} $. For all $u$ belonging to the open set
$\mathscr{O}$ we can
define the quantities $g<u<d$ as the left and right endpoints of the
excursion interval of $\mathscr{O}$ that contains $u$. We define the
\textit{slope} of $c$ at $u$ as the quantity
\[
m_u=\frac{c ( d) -c ( g) }{d-g}=c' ( u) .
\]
The notations $g_u ( f) ,d_u ( f) $ and $m_u ( f) $ will be
preferred when the function~$f$ or the point $u$ are not clear from
context. We will first be interested in continuity properties of the
quantities $g_u$, $d_u$ and $m_u$ when varying the function $f$.

Recall that a~sequence $f_n$ in the space of c\`{a}dl\`{a}g\ functions on
$[0,1]$ converges to $f$ in the Skorohod $J_1$ topology if there exist
a sequence of increasing homeomorphisms from $[0,1]$ into itself such
that $f_n-f\circ\lambda_n$ converges uniformly to $0$ on $[0,1]$.
\begin{lemma}
\label{OneSegmentLemma}
If:
\begin{longlist}[1.]
\item[1.]$f$ is continuous at $u$;
\item[2.]$u\in\mathscr{O}$;
\item[3.] the function
\[
f^l ( t) -\frac{d-t}{d-g}f^l ( g) +\frac{t-g}{d-g}f^l ( d) \qquad\mbox{for $t\in[0,1]$}
\]
is zero only on $\{ g,d\} $;
\item[4.]$f_n\to f$ in the Skorohod $J_1$ topology and $u_n\to u$,
then
\[
g_{u_n} ( f_n) \to g_{u} ( f) ,\qquad d_{u_n} ( f_n) \to d_{u} ( f)
\quad\mbox{and}\quad m_{u_n} ( f_n) \to m_u ( f) .
\]
\end{longlist}
\end{lemma}

The proof is presented in Section~\ref{PathTransformationGeneralCaseSubsection}. We now pass to the analysis
of the particular cases when our L\'{e}vy process $X$ has continuous
paths, or when it reaches its convex minorant continuously.

\subsection{Brownian motion with drift}

In this subsection, we will prove Theorem~\ref{InvarianceUnderPathTransformation} when $X$ is a~(nondeterministic)
L\'{e}vy process with continuous paths, that is, a~(nonzero multiple
of) Brownian motion with drift.

Let $f$ be a~continuous function on $[0,1]$, and consider the
continuous function $\varphi_uf$ given by
%
%
\begin{equation}
\label{TransformationForContinuousFunctions}
\qquad\varphi_uf ( t) =
\cases{
f ( u+t) -f ( u), \hspace*{62pt}\qquad$\mbox{if $0\leq d-u$},$\vspace*{2pt}\cr
f ( d) -f ( u) +f \bigl( g+t-(d-u)\bigr) -f ( g),\vspace*{2pt}\cr \hspace*{140pt}\qquad$\mbox{if $d-u\leq
t\leq d-g$},$\vspace*{2pt}\cr
f ( d) -f ( g) +f \bigl( t-(d-g)\bigr), \qquad$\mbox{if $d-g\leq t\leq
d$},$\vspace*{2pt}\cr
f ( t), \hspace*{117pt}\qquad$\mbox{if $t\leq d$}.$}
\end{equation}
If $f$, $f_n$, $u$ and $u_n$ satisfy the hypotheses of Lemma~\ref{OneSegmentLemma} (which implies that $f_n\to f$ uniformly), then $g (
f_n) \to g ( f) $ and $d ( f_n) \to d ( f) $. Therefore, it
is simple to verify that $( u,f) \mapsto( d-g,\varphi_uf) $ is
continuous at $( u,f) $ when the space of continuous functions on
$[0,1]$ is equipped with the uniform norm. When $X$ is a~L\'{e}vy
process with continuous paths and distributions, that is, a~Brownian
motion with drift, consider its polygonal approximation with step $1/n$
obtained by setting
\[
X^n_{k/n}=X_{k/n} \qquad \mbox{for $k\in\{ 0,1,\ldots, n\} $}
\]
and extending this definition by linear interpolation on $[0,1]$. Then
$X^n\to X$ uniformly on $[0,1]$; it is at this point that the
continuity of the paths of $X$ is important. Now, if $U$ is uniform on
$[0,1]$ and independent of $X$, and we set $U_n=n\lceil U/n \rceil$, then
$( d_n-g_n,\varphi_{U_n}X^n) \to( d-g,\varphi_UX) $. However,
Theorem~\ref{2010AbramsonPitmanTheorem} says that $( U_n,X^n) $ and
$( d_n-g_n,\varphi_{U_n}X^n) $ have the same law. We conclude that
$( U,X) $ and $( d-g,\varphi_U X) $ have the same law, which is
the conclusion of Theorem~\ref{InvarianceUnderPathTransformation} in
this case.

\subsection{Absence of jumps at the convex minorant}

In this subsection, we will prove Theorem~\ref{InvarianceUnderPathTransformation}
when $X$ is a~L\'{e}vy process of
unbounded variation on compact sets [which automatically satisfies
(CD)]. We now let $f$ be a~c\`{a}dl\`{a}g\ function on\vadjust{\goodbreak}
$[0,1]$ and let
$c$ stand for its convex minorant. We will suppose that $f$ is
continuous on the set $\{ c=f^l\} $, which holds whenever $f$ is the
typical trajectory of $X$, thanks to \hyperref[OneJumpPerExcursionProperty]{2}
of Proposition~\ref{BasicPropertiesOfConvexMinorant}.

Again, for all $u\in\{ c<f\} =\{ c<f\wedge f_-\} =\mathscr{O}$ we
define $g$
and $d$ as the left and right endpoints of the excursion interval that
contains $u$.
Since $f$ has jumps, its polygonal approximation does not converge to
it in Skorohod space, but if we define
\[
f_n ( t) =f ( [nt]/n) ,
\]
then $f_n$ converges in the Skorohod $J_1$ topology to $f$ as $n\to
\infty$; cf. \citeauthor{MR1700749} (\citeyear{MR1700749}), Chapter 2, Lemma 3, page 127.
This will
called the piecewise constant approximation to $f$ with span $1/n$ and
is the way we will choose to approximate a~L\'{e}vy process when it has
jumps. The first complication in this case is that the discrete
invariant path transformation was defined for the polygonal
approximation and not for the piecewise constant approximation to our
L\'{e}vy process. For this reason, we will have to define a~more
flexible path transformation than in the continuous case: for
$u_1<u_2<u_3\in(0,1)$, we define $\varphi_{u_1,u_2,u_3}f$ by
%
%
\begin{equation}
\label{TransformationForCadlagFunctionsNotJumpingAtCM}
\qquad\varphi_{u_1,u_2,u_3} f ( t) =
\cases{
f ( u_2+t) -f ( u_2), \qquad 0\leq t< u_3-u_2,\vspace*{2pt}\cr
f ( u_3) -f ( u_2)+f \bigl( u_1+t-(u_3-u_2)\bigr) -f ( u_1),
\vspace*{2pt}\cr
\hspace*{110pt} u_3-u_2\leq t\leq
u_3-u_1,\vspace*{2pt}\cr
f ( u_3) -f ( u_1)+f \bigl( t-(u_3-u_1)\bigr), \vspace*{2pt}\cr
\hspace*{110pt}u_3-u_1\leq
t<u_3,\vspace*{2pt}\cr
f ( t), \hspace*{63pt}\qquad u_3\leq t.}
\end{equation}
The path transformation $\varphi_u$ of (\ref{TransformationForContinuousFunctions}) corresponds to $\varphi
_{g,u,d}$. We are interested in $\varphi_{g,U,d}X$, which will be
approximated by $\varphi_{\tilde g_n,U_n,\tilde d_n}X^n$ where $\tilde
g_n$ and $\tilde d_n$\vspace*{1pt} are the left and right endpoints of the excursion
of the polygonal approximation to $X$ of span $1/n$ which contains
$U_n=\lceil Un \rceil/n$, and $X^n$ is the piecewise constant
approximation to
$X$ with span $1/n$. We are forced to use both the vertices of the
convex minorant of the polygonal approximation and the piecewise
constant approximation, since $X^n\to X$ (in the Skorohod $J_1$
topology), but with $( \tilde g_n,\tilde d_n) $ we can define a~nice invariant transformation: Theorem~\ref{2010AbramsonPitmanTheorem}
asserts that
\[
( U_n,\varphi_{\tilde g_n,U_n,\tilde d_n}X^n) \quad\mbox{and}\quad( \tilde
d_n-\tilde g_n,X^n)
\]
have the same law.
Indeed, Theorem~\ref{2010AbramsonPitmanTheorem}
is an assertion about the increments of a~random walk and the polygonal
and piecewise approximations to $X$ of span $1/n$ are constructed from
the same increments. 

Lemma~\ref{OneSegmentLemma} tells us that $( \tilde g_n,\tilde d_n)
\to( d,g) $. It is therefore no surprise that
\[
\varphi_{\tilde g_n,U_n,\tilde d_n}X^n\to\varphi_{g,U,d}X,
\]
telling us that $( U,X) $ and $( d-g,\varphi_{g,U,d}X) $ have
the same law whenever $X$ satisfies~(CD) and has unbounded
variation on finite intervals. Convergence follows from the following
continuity assertion:\vadjust{\goodbreak}
\begin{lemma}
\label{SkorohodContinuityLemmaForPathTransformationWithoutJumps}
If $f$ is continuous at $( u_1,u_2,u_3) $, $f_n\to f$ in the
Skorohod $J_1$ topology, and $u^n_i\to u_i$ for $i=1,2,3$, then
\[
\varphi_{u^n_1,u^n_2,u^n_3}f_n\to\varphi_{u_1,u_2,u_3}f.
\]
\end{lemma}

Lemma~\ref{SkorohodContinuityLemmaForPathTransformationWithoutJumps} is
an immediate consequence of the following convergence criterion found
in \citeauthor{MR838085} (\citeyear{MR838085}), Proposition III.6.5, page 125.
\begin{proposition}
\label{EthierKurtzCriterionForWeakConvergence}
A sequence $f_n$ of c\`{a}dl\`{a}g\ functions on $[0,1]$ converges to
$f$ in
the Skorohod $J_1$ topology if and only if for every sequence $( t_n)
\subset[0,1]$ converging to $t$:
\begin{longlist}[1.]
\item[1.]\label{EKCondition1} $|f_n ( t_n) -f ( t) |\wedge|f_n ( t_n) -f
( t-) |\to0$;
\item[2.]\label{EKCondition2} if $|f_n ( t_n) -f ( t) |\to0$,
$t_n\leq s_n\to t$, then $|f_n ( s_n) -f ( t) |\to0$;
\item[3.]\label{EKCondition3} if $|f_n ( t_n) -f ( t-) |\to0$,
$s_n\leq t_n$ and $s_n\to t$, then $|f_n ( s_n) -f ( t) |\to0$.
\end{longlist}
\end{proposition}

In particular, we see that if $f$ is continuous at $t$, then $f_n (
t_n) \to f ( t) $. The above criterion is clearly necessary for
convergence since if $f_n\to f$, then there exist a~sequence $( \lambda
_n,n\in\mathbb{N}) $ of increasing homeomorphisms of $[0,1]$ into
itself such that
$f_n-f\circ\lambda_n$
converges to zero uniformly. If $t_n\to t$, then $f_n ( t_n) $ will
be close to either $f ( t-) $ or $f ( t) $ depending on if $\lambda_n
( t_n) <t$ or $\lambda_n ( t_n) \geq t$. By using the above
criterion, we focus on the real problem for continuity for the
transformation $\varphi_{u_1,u_2,u_3}$, namely, that nothing wrong
happens at $u_3-u_2$, $u_3-u_1$ and $u_3$.
\begin{pf*}{Proof of Lemma
\protect\ref{SkorohodContinuityLemmaForPathTransformationWithoutJumps}}
Let us prove that for every $t\in[0,1]$, the conditions of Proposition
\ref{EthierKurtzCriterionForWeakConvergence} hold for $\varphi
_{u^n_1,u^n_2,u^n_3}f_n$ and $\varphi_{u_1,u_2,u_3}f$.\vspace*{1pt}

Let $\lambda_n$ be increasing homeomorphisms of $[0,1]$ into itself
such that
\[
f_n-f\circ\lambda_n\to0
\]
uniformly. We proceed by cases.
\begin{longlist}
\item[$t<u_3$.] Eventually $t<u^n_3$, so that $\varphi
_{u^n_1,u^n_2,u^n_3}f_n ( t) =f_n ( t) $ and $\varphi _{u_1,u_2,u_3}f
( t) =f ( t) $. Since $f_n$ and $f$ satisfy the
conditions of Proposition~\ref{EthierKurtzCriterionForWeakConvergence}
at time $t$, the same holds for their images under the path transformation.
\item[$t<u_3-u_2$.] Eventually $t<u^n_3-u^n_2$ so that
\[
\varphi_{u_1,u_2,u_3}f ( t) =f ( u_2+t) -f ( u_2)\quad \mbox{and} \quad\varphi_{u^n_1,u^n_2,u^n_3}f^n ( t) =f^n ( u^n_2+t) -f^n (
u^n_2) .
\]
Since $f$ is continuous at $u_2$, Proposition~\ref{EthierKurtzCriterionForWeakConvergence} implies that $f^n ( u^n_2) \to
f ( u) $, so that
$\varphi_{u^n_1,u^n_2,u^n_3}f^n ( t)$ can be made arbitrarily
close to either $\varphi_{u_1,u_2,u_3}f ( t)$ or $\varphi
_{u_1,u_2,u_3}f ( t-)$
depending on if
\[
u_n+t_n<\lambda_n^{-1} ( u+t)\quad \mbox{or} \quad u_n+t_n\geq\lambda_n^{-1} (
u+t) .
\]
\item[$t\in(u_3-u_2,u_3)\setminus\{ u_3-u_1\} $] is analogous to the
preceding case.
For
$t\in\{ u_3-u_2,u_3-u_1,u_3\}$ set
\[
v_1=u_3-u_2,\qquad v_2=u_3-u_1\quad\mbox{and}\quad v_3=u_3.\vadjust{\goodbreak}
\]
Since $f$ is continuous at $u_3$, condition \hyperref[EKCondition3]{3} gives
\[
f_n ( u_i^n) \to f ( u_i)\qquad \mbox{for $i=1,2,3$,}
\]
and so
\[
\varphi_{u_1,u_2,u_3}f ( v^n_i) \to\varphi _{u_1,u_2,u_3}f ( v_i)
\qquad\mbox{for $i=1,2,3$}.
\]
\end{longlist}
\upqed\end{pf*}

\subsection{The general case}\label{PathTransformationGeneralCaseSubsection}

In this subsection, we prove Theorem~\ref{InvarianceUnderPathTransformation} for a~L\'{e}vy process $X$ under
the sole assumption (CD).

The challenge to overcome in the remaining case, in which $X$ can jump
into and out of the convex minorant, is to show how one can handle the
jumps; although a~result in the vein of Lemma~\ref{SkorohodContinuityLemmaForPathTransformationWithoutJumps} will play a~prominent role in our analysis, a~more careful inspection of how $g_n$
differs from $g$ is needed in order to sort the following problem: in
general, the operation of rearranging pieces of c\`{a}dl\`{a}g\ paths
is not
continuous and depends sensitively on the points at which the
rearrangement is made. A simple example helps to clarify this: consider
$f={\mathbf{1}}_{[1/3,1]}+{\mathbf{1}}_{[2/3,1]}$, so that if
$u_1=1/3, u_2=1/2$ and
$u_3=2/3$, we have $\varphi_{u_1,u_2,u_3}f={\mathbf
{1}}_{[1/6,1]}+{\mathbf{1}}_{[2/3,1]}$. Note that if $u^n_1\to u_1$
and $u^n_1\in(0,1/2)$, then
\[
\varphi_{u^n_1,u_2,u_3}f=
\cases{
{\mathbf{1}}_{[1/6,1]}+{\mathbf{1}}_{[1/6+1/3-u^n_1,1]},&\quad $\mbox{if
$u^n_1\in(0,1/3]$},$\vspace*{2pt}\cr
{\mathbf{1}}_{[1/6,1]}+{\mathbf{1}}_{[1-u^n_1,1]},&\quad $\mbox{if $u^n_1\in
[1/3,1/2)$}.$}
\]
We conclude that $\varphi_{u^n_1,u_2,u_3}f\to\varphi_{u_1,u_2,u_3}f$ if
and only if $u^n_1\geq1/3$ eventually.

Let $f\dvtx [0,1]\to\mathbb{R}$ be a~c\`{a}dl\`{a}g\ function which starts
at zero
and $c$ its convex minorant on $[0,1]$.
Let also $f^l=f\wedge f_-$ be the lower semicontinuous regularization
of $f$. As before, the component intervals of the open set $\mathscr
{O}=\{ c<f^l\} $ are called the excursion intervals of $f$, and that
for $u\in
\mathscr{O}$,
$( g,d) $ is the excursion interval that contains $u$.

\begin{figure}

\includegraphics{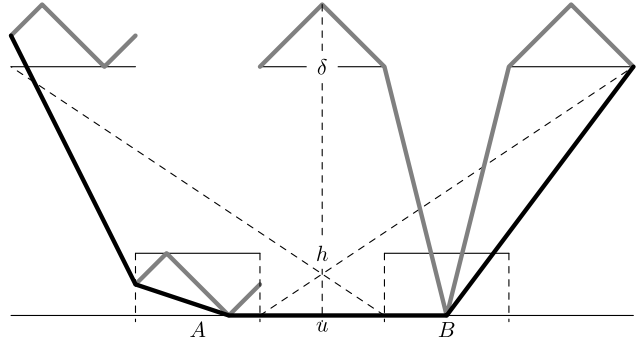}

\caption{Visual content of Lemma \protect\ref{VisualCoreOfOneSegmentLemma}.}
\label{VisualCoreOfOneSegmentLemmaFigure}\vspace*{3pt}
\end{figure}

We first give the proof of Lemma~\ref{OneSegmentLemma}; the proof
depends on another lemma with a~visual appeal, which is to be
complemented with Figure~\ref{VisualCoreOfOneSegmentLemmaFigure}.
%
%

\begin{lemma}
\label{VisualCoreOfOneSegmentLemma}
If for a~c\`{a}dl\`{a}g\ function $f\dvtx [0,1]\to\mathbb{R}$:
\begin{longlist}[1.]
\item[1.] there exist closed intervals $A$ and $B$ in $[0,1]$ such that
$\inf B-\sup A>0$ and
\item[2.] there exists $\delta>0$ and
\[
h<\delta\frac{\inf B-\sup A}{\inf B\vee( 1-\sup A) }
\]
such that
\[
f>\delta\mbox{ on }[0,1]\setminus A\cup B \quad\mbox{and}\quad
\min_{x\in A\cup B} f^l ( x) <h,
\]
\end{longlist}
then for all $u\in(\sup A,\inf B)$,
\[
g_u\in A,\qquad d_u\in B \quad\mbox{and}\quad m_u\leq\frac{h}{\inf
B-\sup A}.\vadjust{\goodbreak}
\]
\end{lemma}
\begin{pf}
This assertion can be checked by cases. We consider $3$ possible
positions for $g_u$ and three other for $d_u\dvtx  g_u<\inf A$, $g_u\in A$
and $g_u\in(\sup A,u)$ and similarly $d_u\in(u,\inf B)$, $d_u\in B$
and $d_u>\sup B$. We number each from~1 to 3 and write $C_{i,j}$ for
the corresponding case. We trivially discard the cases
\[
C_{1,1}, C_{1,3}, C_{3,1}, C_{3,3}
\]
for each one would force $c ( g) $ to be above the zero slope line
through $( 0,\delta) $, hence to pass above $g$ on $A$ and $B$. The
case $C_{2,1}$ would force $c$ (hence $f$) to be above~$\delta$ on $B$
while $C_{3,2}$ would force $f$ to be above $\delta$ on $A$, hence both
are discarded. We finally discard the case $C_{2,3}$ (and by a~similar
argument~$C_{3,2}$) because of our choice of $h$, since a~line from a~point of $A\times[0,h]$ to $[\sup B,1]\times[\delta,\infty)$ passes
above $h$ on $B$.
\end{pf}

\begin{pf*}{Proof of Lemma \protect\ref{OneSegmentLemma}}
Set $u\in\{ c<f^l\} $, and write $g$ and $d$ for $g_u ( f) $ and
$d_u ( f) $ so that $g<u<d$. Recall that $c$ is linear on $(g_u ( f)
,d_u ( f) )$.
By considering instead
\begin{eqnarray*}
t&\mapsto& f ( t) -\frac{d-t}{d-g}f^l ( g) +\frac{t-g}{d-g}f^l ( d)
\quad\mbox{and}\\
 t&\mapsto& f_n ( t) -\frac{d-t}{d-g}f^l ( g) +\frac
{t-g}{d-g}f^l ( d) ,
\end{eqnarray*}
our assumptions allow us to reduce to the case
\[
f^l ( g) =f^l ( d) =0\quad \mbox{and} \quad f^l>0\qquad\mbox{on
$[0,1]\setminus\{ g,d\} $}.
\]
We will now consider the case $0<g<d<1$, the cases $g=0$ or $d=1$ being
handled similarly.

For every
\[
\varepsilon<g\wedge(1-d)\wedge\frac{d-g}{2}
\]
we can define
\begin{eqnarray*}
\delta ( \varepsilon)
&=&\inf\{ f ( t) \dvtx  t\in[0,g-\varepsilon]\cup[g+\varepsilon
,d-\varepsilon]\cup[d+\varepsilon,1]\}
\\
&=&\min\{ f^l ( t) \dvtx  t\in[0,g-\varepsilon]\cup[g+\varepsilon
,d-\varepsilon]\cup[d+\varepsilon,1]\} .
\end{eqnarray*}
Then $\delta ( \varepsilon) >0$ for $\varepsilon>0$ and $\delta (
\varepsilon) \to0$
as $\varepsilon\to0$. Since $g<u<d$, we can choose $\varepsilon$
small enough so that
\[
u\in(g+\varepsilon,d-\varepsilon).
\]
Since $f_n\to f$, there exists a~sequence of increasing homeomorphisms
$\lambda_n$ of $[0,1]$ converging uniformly to the identity such that
\[
f_n-f\circ\lambda_n
\]
converges uniformly to zero. (If $f$ is continuous, $\lambda_n$ can be
taken equal to the identity function.)

Also, given $h_n$ eventually bounded away from $0$,
\[
\min\{ f^l ( t) \dvtx t\in(g-\varepsilon,g+\varepsilon)\} <h_n
\quad\mbox{and}\quad
\min\{ f^l ( t) \dvtx t\in(d-\varepsilon,d+\varepsilon)\} <h_n
\]
for large enough $n$. Hence
\[
\min\bigl\{ f^l_n ( t) \dvtx t\in\bigl(\lambda_n^{-1} ( g-\varepsilon) ,\lambda
_n^{-1} ( g+\varepsilon) \bigr)\bigr\} <h_n
\]
%
and
\[
\min\bigl\{ f^l_n ( t) \dvtx t\in\bigl(\lambda_n^{-1} ( d-\varepsilon) ,\lambda
_n^{-1} ( d+\varepsilon) \bigr)\bigr\} <h_n
\]
for large enough $n$. The particular $h_n$ we will consider is
\begin{eqnarray*}
h_n&=& \delta ( \varepsilon) \frac{\lambda_n^{-1} ( d-\varepsilon)
-\lambda _n^{-1} ( g+\varepsilon) }{\lambda_n^{-1} ( d-\varepsilon)
\vee(1-\lambda _n^{-1} ( g+\varepsilon) )}
\\
& \to&\delta ( \varepsilon) \frac{( d-g-2\varepsilon) }{( (
d-\varepsilon ) \vee(1-g-\varepsilon)) }>0
\end{eqnarray*}
which is eventually positive. Since $f>\delta ( \varepsilon) $ on
$[0,g-\varepsilon
]\cup[g+\varepsilon,d-\varepsilon]\cup\allowbreak[d+\varepsilon,1]$, then
\[
f_n>\delta\qquad \mbox{on } [0,\lambda_n^{-1} ( g-\varepsilon) ]\cup
[\lambda_n^{-1} ( g+\varepsilon) ,\lambda_n^{-1} ( d-\varepsilon)
]\cup[\lambda _n^{-1} ( d+\varepsilon) ,1],
\]
and Lemma~\ref{VisualCoreOfOneSegmentLemma} now tells us that
\begin{eqnarray*}
g_{u_n} ( f_n) &\in&\bigl(\lambda_n^{-1} ( g-\varepsilon) ,\lambda
_n^{-1} ( g+\varepsilon) \bigr),\\
d_{u_n} ( f_n) &\in&\bigl(\lambda_n^{-1} ( d-\varepsilon) ,\lambda
_n^{-1} ( d+\varepsilon) \bigr)
\end{eqnarray*}
and
\[
m_{u_n} ( g_n) \leq h_n/\bigl(\lambda_n ( d-\varepsilon) -\lambda _n (
g+\varepsilon) \bigr),
\]
so that eventually
\begin{eqnarray*}
g_{u_n} ( f_n) &\in&(g-2\varepsilon,g+2\varepsilon),\qquad d_{u_n} ( f_n)
\in
(d-2\varepsilon,d+2\varepsilon)\quad \mbox{and}\\
 m_{u_n} ( f_n) &\leq&
2\delta/
(d-\varepsilon)\vee(1-g-\varepsilon).
\end{eqnarray*}
\upqed\end{pf*}

\begin{remark*}
In the context of the above proof, if we suppose that
$f ( g-) =c ( g) <f ( g) $ and $f ( d-) =c ( d) $, then
for $h_n$ eventually bounded away from zero, we actually have
\[
\min\{ f^l ( t) \dvtx t\in[g-\varepsilon,g)\} <h_n
\]
for large enough $n$, and so we get
\[
g_{u_n} ( f_n) <\lambda_n^{-1} ( g) .
\]
%
This remark is crucial to the proof of Theorem~\ref{InvarianceUnderPathTransformation}.
\end{remark*}
\begin{remark*}
Let $c_n$ be the convex minorant of $f_n$. Under the hypotheses of
Lemma~\ref{OneSegmentLemma}, we can actually deduce that if $t_n\to g$,
then $c_n ( t_n) \to c ( g) $, while if $t_n\to d$, then $c_n ( t_n)
\to c ( d) $. This is because of the following result about
convergence of convex functions.
\end{remark*}
\begin{proposition}
If $c_n$ and $c$ are convex functions on $[0,1]$, for some $a\in(0,1)$
we have $c_n ( a) \to c ( a) $, and if the two sequences $( c_n ( 0) )
$ and $( c_n ( 1) ) $ are bounded, then for every
sequence $a_n\to a$ we have $c_n ( a_n) \to c ( a) $.
\end{proposition}
\begin{pf}
If $a_n\leq a$, we can use the inequalities
\[
c_n ( a_n) \leq c_n ( a) \frac{a_n}{a}+c_n ( 0) \frac{a-a_n}{a}
\]
and
\[
c_n ( a_n) \geq c_n ( a) \frac{1-a}{1-a_n}+c_n ( 1) \frac
{a-a_n}{1-a_n}.
\]
We get an analogous pair of inequalities when $a_n\geq a$, which allows
us to conclude that the sequence $( c_n ( a) -c_n ( a_n) ) $
goes to zero.
\end{pf}

Given $u_1<u_2<u_3$, we now define a~new c\`{a}dl\`{a}g\ function $\psi
_{u_1,u_2,u_3} f$ as follows:
\[
\label{GeneralDefinitionOfPathTransformation}
\psi_{u_1,u_2,u_3} f ( t) =
\cases{
f ( u_2+t) -f ( u_2), \hspace*{66pt}\qquad 0\leq t< u_3-u_2,\vspace*{2pt}\cr
c ( u_3) -c ( u_1) +f \bigl( u_1+t-(u_3-u_2)\bigr)-f ( u_2),\vspace*{2pt}\cr
 \hspace*{154pt}\qquad u_3-u_2\leq t\leq
u_3-u_1,\vspace*{2pt}\cr
c ( u_3) -c ( u_1) +f \bigl( t-(u_3-u_1)\bigr), \qquad u_3-u_1\leq
t<u_3,\vspace*{2pt}\cr
f ( t), \hspace*{130pt}\qquad\mbox{$u_3\leq t$}.}
\]
The difference with the path transformations of (\ref{TransformationForContinuousFunctions}) and (\ref{TransformationForCadlagFunctionsNotJumpingAtCM}) is that we now use the
convex minorant $c$ instead of only the function $f$. This has the
effect of choosing where to place the jumps that $f$ might make as it
approaches its convex minorant. Note, however, that $\psi
_{u_1,u_2,u_3}f=\varphi_{u_1,u_2,u_3}f$ if $f=c$ at~$u_1$ and $u_3$.

Our next task will be to analyze the continuity of $f\mapsto\psi
_{g,u,d}f$ on Skorohod space, with special emphasis on the
approximations we will use. For every~$n$,~$f_n$ and $\tilde f_n$ will
be the piecewise constant and polygonal approximations to~$f$ with span
$1/n$, we set $u_n=\lceil nu \rceil/n$, and
\[
g_n=g_{u_n} ( f_n) , \qquad d_n=d_{u_n} ( f_n) ,\qquad \tilde
g_n=g_{u_n} ( \tilde f_n) \quad \mbox{and}\quad \tilde d_n=d_{u_n} ( \tilde
f_n) .
\]

\begin{lemma}
\label{ContinuityOfPathTransformationLemma}
Under the hypotheses of Lemma~\ref{OneSegmentLemma}, if
either
\[
f ( g) =c ( g) \quad\mbox{and}\quad f ( d) =c ( d)
\quad\mbox{or}\quad
f ( g-) =c ( g) \quad\mbox{and}\quad f ( d-) =c ( d) ,
\]
then
\[
\psi_{\tilde g_n,u_n,\tilde d_n}f_n\to\psi_{g,u,d}f
\]
in the Skorohod $J_1$ topology.
\end{lemma}

\begin{pf}
Since we have already analyzed what happens when $f$ is continuous at
$g$ and $d$, the essence of the argument will be illustrated when
\[
f ( g) =c ( g) <f ( g-) \quad\mbox{and}\quad f ( d) =c ( d) <f ( d-) .
\]
As in the proof of Lemma~\ref{SkorohodContinuityLemmaForPathTransformationWithoutJumps}, we verify
that for every $t\in[0,1]$, the conditions of Proposition~\ref{EthierKurtzCriterionForWeakConvergence} hold for $\psi_{\tilde
g_n,u_n,\tilde d_n}f_n$ and $\psi_{d,u,g}f$ at time $t$.\vspace*{1pt}

Let $\lambda_n$ be a~sequence of increasing homeomorphisms of $[0,1]$
such that
\[
f_n-f\circ\lambda_n\to0
\]
uniformly. The crucial part of the argument is to use the remarks after
Lemma~\ref{OneSegmentLemma} from which we deduce that
\[
\lambda^{-1}_n ( g) \leq g_n \quad\mbox{and}\quad \lambda ^{-1}_n ( d) \leq d_n.
\]
Since $f_n$ is the piecewise constant approximation to $f$, then
$\lambda_n^{-1}$ must eventually take $g$ and $d$ to $[ng]/n$ and
$[nd]/n$. But comparing the convex minorants of the piecewise
constant and polygonal approximations to $f$ with span $1/n$ leads to
\[
g_n-1/n\leq\tilde g_n \quad\mbox{and}\quad d_n-1/n\leq\tilde d_n
\]
so that
\[
\lambda_n^{-1} ( g) \leq\tilde g_n\quad \mbox{and}\quad \lambda _n^{-1} ( d)
\leq\tilde d_n.
\]
Again using the remarks after the proof of Proposition~\ref{TransformationForCadlagFunctionsNotJumpingAtCM}, we see that
\[
c_n ( \tilde g_n) \to c ( g) \quad\mbox{and}\quad c_n ( \tilde d_n) \to c ( d) .
\]

The conditions of Proposition~\ref{EthierKurtzCriterionForWeakConvergence} can now be verified at times
$t\in[0,1]\setminus{d-u,d}$ as in the proof of Lemma~\ref{TransformationForCadlagFunctionsNotJumpingAtCM}, while for $t\in\{
d-u,d\} $, the proof is similar and hence will be illustrated when
$t=d-g$. Since $f_n-f\circ\lambda_n\to0$ uniformly,\vadjust{\goodbreak} the jump of~$f$ at
$g$ is approximated by the jump of $f_n$ at $\lambda_n^{-1} ( g) $.
We reduce to cases by taking subsequences: when $t_n>\tilde d_n-u_n$
for all $n$, then $t_n+u_n>\lambda_n^{-1} ( d) $ so that
\[
\psi_{\tilde g_n,u_n,\tilde d_n}f ( t_n) \to f ( d) -f ( g) +f ( g) -f
( u) =f ( d) -f ( u) .
\]
On the other hand, when $t_n\leq d_n-u_n$ for all $n$, we see that
\[
\psi_{\tilde g_n,u_n,\tilde d_n}f ( t_n) \mbox{ is close to }
f ( d-) -f ( u)
\mbox{ or }
f ( d) -f ( u)
\]
depending on if
\[
t_n+u_n< \lambda_n^{-1} ( d)
\quad\mbox{or}\quad
t_n+u_n\geq\lambda_n^{-1} ( d) .
\]
Hence, the conditions of Proposition~\ref{EthierKurtzCriterionForWeakConvergence} are satisified at $t=d-u$.
\end{pf}


We finally proceed to the proof of Theorem~\ref{InvarianceUnderPathTransformation}.
\begin{pf*}{Proof of Theorem
\protect\ref{InvarianceUnderPathTransformation}}
Thanks to Proposition~\ref{BasicPropertiesOfConvexMinorant}, $X$ almost
surely satisfies the conditions of Lemma~\ref{ContinuityOfPathTransformationLemma} at $U$. Hence, $( d_n-g_n,\psi
_{\tilde g_n,U_n,\tilde d_n} ( X^n) ) $ converges in law
to $( d-g,\psi_{d,U,g}X) $ thanks to Lemmas~\ref{OneSegmentLemma}
and~\ref{ContinuityOfPathTransformationLemma}, as well as the
continuous mapping theorem. Since $( U_n,X^n) $ converges in law to
$( U,X) $ and the laws of $( U_n,X^n) $ and $( d_n-g_n,\psi_{\tilde
g_n,U_n,\tilde d_n} ( X_n) ) $ are equal by Theorem~\ref{2010AbramsonPitmanTheorem}, then $( U,X) $ and $( d-g,X^U) $
have the same law.
\end{pf*}

\section{Excursions above the convex minorant on a~fixed interval}\label{ExcursionSection}

In this section we will prove Theorem~\ref{PointProcessOfExcursions},
which states the equality in law between two sequences. We recall the
setting: $X$ is a~L\'{e}vy process such that $X_t$ has a~continuous
distribution for every $t>0$, $C$ is its convex minorant on $[0,1]$,
$X^l=X\wedge X_-$ is the lower semicontinuous regularization of $X$,
$\mathscr{O}=\{ C<X^l\} $ is the excursion set, $\mathscr{I}$ is the
set of
excursion intervals of $\mathscr{O}$, for each $( g,d) \in\mathscr
{I}$, and
we let $e^{( g,d) }$ be the excursion associated to $( g,d) $
given by
\[
e^{( g,d) }_s=X_{( g+s) \wedge d}-C_{( g+s) \wedge d}.
\]
We ordered the excursion intervals to state Theorem~\ref{PointProcessDescriptionForDeterministicAndFiniteHorizon} by sampling
them with an independent sequence of uniform random variables on $[0,t]$.

The first sequence of interest is
\[
\bigl( \bigl( d_i-g_i,C_{d_i}-C_{g_i},e^{( g_i,d_i) }\bigr) ,i\geq1\bigr) .
\]
The second sequence is obtained with the aid of an independent
stick-breaking process and the Vervaat transformation. Recall that
$V_tf$ stands for the Vervaat transform of $f$ on $[0,t]$. Let
$V_1,V_2,\ldots$ be an i.i.d. sequence of uniform random variables on
$(0,1)$, and construct
\[
L_1=V_1, \qquad L_n=V_n( 1-V_1) \cdots( 1-V_{n-1}) \quad\mbox
{and}\quad S_i=L_1+\cdots+L_i.
\]
This sequence helps us to break up the paths of $X$ into the
independent pieces $Y^i$, $i=1,2,\ldots$ given by
\[
Y^i_t=X_{S_{i-1}+t}-X_{S_{i-1}},\qquad 0\leq t\leq L_i,\vadjust{\goodbreak}
\]
from which we can define the sequence of Knight bridges,
\[
K^i_t=Y^i_t-\frac{t}{L_i}Y^i_{L_i},\qquad 0\leq t\leq L_i.
\]
Our second sequence is
\[
\bigl( \bigl( L_i,X_{S_i}-X_{S_{i-1}},V_{L_i} ( K^i) \bigr) ,i\geq1\bigr) .
\]

To prove the equality in law, we will use Theorem~\ref{InvarianceUnderPathTransformation}
to obtain a~process~$\tilde X$
which has the same law as $X$, as well as a~stick-breaking sequence
$\tilde L$ independent of $\tilde X$ such that, with analogous
notation, the pointwise equality
\[
\bigl( \bigl( d_i-g_i,C_{d_i}-C_{g_i},e^{( g_i,d_i) }\bigr) ,i\geq1\bigr) =
\bigl( \bigl( L_i,\tilde X_{ \tilde S_i}-\tilde X_{\tilde S_{i-1}},V_{\tilde L_i}
( \tilde K^i) \bigr) ,i\geq1\bigr)
\]
holds. This proves Theorem~\ref{PointProcessOfExcursions}.

Let us start with the construction of $\tilde X$ and $\tilde L$. Apart
from our original L\'{e}vy process $X$, consider an i.i.d. sequence of
uniform random variables $U_1,U_2,\ldots$ independent of $X$.
Consider first the connected component $( g_1,d_1) $ of $\{ C<X\wedge
X_-\} $ which contains $U_1$ and let $X^1$ be the result of
applying the path transformation of Theorem~\ref{InvarianceUnderPathTransformation} to $X$ at the points $g_1$, $U_1$
and $d_1$. We have then seen that $\tilde V_1=d_1-g_1$ is uniform on
$[0,1]$ and independent of $X^1$. Set $\tilde S_0=0$ and $\tilde
L_1=\tilde V_1$.

Consider now the convex minorant $C^1$ of
\[
Z^1=X^1_{\tilde L_1+\cdot}-X^1_{\tilde L_1}
\]
on $[0,1-\tilde L_1]$: we assert that it is obtainable from the graph
of $C$ by erasing the interval $( g_1,d_1) $ and closing up the
gap, arranging for continuity. Formally, we assert the equality
\[
C^1_t=
\cases{
C_{t},&\quad$\mbox{if }t\in[0,g_1),$\vspace*{2pt}\cr
C_{t-g_1+d_1}-( C_{d_1}-C_{g_1}),&\quad $\mbox{if }t\in[g_1,1-\tilde L_1].$}
\]
%
Note that $C^1$ is continuous on $[0,1-\tilde L_1]$ by construction and
it is convex by a~simple analysis.
To see that $C^1$ is the convex minorant of $Z^1$, we only need to
prove that at $g_1$ it coincides with $Z^1_{g_1}\wedge Z^1_{g_1-}$; cf.
Figure~\ref{PathTransformation2} to see how it might go wrong. 
If $X_{d_1}=C_{d_1}$, then
\[
Z^1_{g_1}
=X_{d_1}-( C_{d_1}-C_{g_1}) =C_{g_1}=C^1_{d_1},
\]
while if $X_{d_1-}=C ( d_1) <X_{d_1}$, then property \hyperref[OneJumpPerExcursionProperty]{2} of Proposition~\ref{BasicPropertiesOfConvexMinorant} implies that $X_{g_1-}=C ( g_1) $ and
\[
Z^1_{d_1-}=X^1_{d_1-}-X^1_{\tilde L_1}=C_{d_1}-C_{g_1}+X_{g_1-}-(
C_{d_1}-C_{g_1}) =C_{g_1}=C^1_{g_1}.
\]


Let $( g_2,d_2) $ be the connected component of $\{ C^1<Z^1\} \subset
[0,1-\tilde L_1]$ that contains $U_2( 1-\tilde L_1) $ and define
\[
\tilde L_2=d_2-g_2,\qquad \tilde V_2=\frac{d_2-g_2}{1-\tilde V_1}
\]
as well as the process $X^2$ which will be the concatenation of $X^1$
on $[0,\tilde V_1]$ as well as the path transformation of $Z^1$ on
$[0,1-\tilde V_1]$; that is, $Z^1$ transformed according to the path
transformation of Theorem~\ref{InvarianceUnderPathTransformation} with
parameters~$g_2,U_2( 1-\tilde L_1) ,d_2$. From Theorem~\ref{InvarianceUnderPathTransformation} and the independence of
\[
Z^2=X^1_{\cdot+\tilde V_1}-X^1_{\tilde V_1} \quad\mbox{and}\quad
X^1_{\cdot\wedge\tilde V_1}
\]
we see that:
\begin{longlist}[1.]
\item[1.]$X^2$ has the same law as $X^1$;
\item[2.]$\tilde V_1$ and $\tilde V_2$ are independent of $X^2$, and
$\tilde V_2$ is independent of $\tilde V_1$ and has an uniform
distribution on $(0,1)$;
\item[3.] the convex minorant $C^2$ of $Z^2$ on $[0,1-\tilde L_1-\tilde
L_2]$ is obtained from $C^1$ by deleting the interval $(g_2,d_2)$ and
closing up the gap arranging for continuity.
\end{longlist}
Now it is clear how to continue the recursive procedure to obtain, at
step $n$ a~sequence $\tilde V_1,\ldots,\tilde V_n$ and a~process $X^n$
such that if $\tilde L_n=\tilde V_n( 1-\tilde V_{n-1}) \cdots(
1-\tilde V_1) $ and $\tilde S_n=\tilde L_1+\cdots+\tilde L_n$; then:
\begin{longlist}
\item[1.]$X^n$ has the same law as $X$.
\item[2.]$X^n$, $\tilde V_1,\ldots,\tilde V_n$ are independent an the
latter $n$ variables are uniform on $(0,1)$.
\item[3.] Let $C^n$ is the convex minorant of
\[
Z^n=X^n_{\tilde S_n+\cdot}-X^n_{\tilde S_n}
\]
on $[0,1-\tilde S_n]$. Then $C^n$ is obtained from $C^{n-1}$ by
removing the selected interval $( g_n,d_n) $ and closing up the gap
arranging for continuity.
%
\item[4.]\label{ProjectivePropertyOfIteratedPathTransformation} $X^n$
coincides with $X^{n-1}$ on $[0, \tilde S_{n-1}]$.
\end{longlist}

From property \hyperref[ProjectivePropertyOfIteratedPathTransformation]{4}
above, it is clear that $X^n$ converges pointwise on $[0,1]$ almost
surely: it clearly does on $[0,1)$ and $X^n_1=X_1$. Also, we see that
$\tilde X$ has the same law as $X$ and that it is independent of
$V_1,V_2,\ldots,$ which is an i.i.d. sequence of uniform random variables.

\section*{Acknowledgments}
Josh Abramson and Nathan Ross provided valuable input to the present
work, not only through enlightening discussions, but also by providing
access to their recent work on convex minorants.

We gladly thank the anonymous referee for his expeditious revision.

\bibliographystyle{imsart-nameyear}

\begin{thebibliography}{44}

\bibitem[\protect\citeauthoryear{Abramson and
  Pitman}{2011}]{2010AbramsonPitman}
\begin{barticle}[author]
\bauthor{\bsnm{Abramson},~\bfnm{Josh}\binits{J.}} \AND
  \bauthor{\bsnm{Pitman},~\bfnm{Jim}\binits{J.}}
(\byear{2011}).
\btitle{Concave majorants of random walks and related Poisson processes}.
\bjournal{Combinatorics Probab. Comput.}
\bvolume{20}
\bpages{651--682}.
\end{barticle}
\endbibitem

\bibitem[\protect\citeauthoryear{Abramson et~al.}{2011}]{Abramson2011fk}
\begin{barticle}[author]
\bauthor{\bsnm{Abramson},~\bfnm{Josh}\binits{J.}},
  \bauthor{\bsnm{Pitman},~\bfnm{Jim}\binits{J.}},
  \bauthor{\bsnm{Ross},~\bfnm{Nathan}\binits{N.}} \AND
  \bauthor{\bsnm{Uribe~Bravo},~\bfnm{Ger{\'{o}}nimo}\binits{G.}}
(\byear{2011}).
\btitle{Convex minorants of random walks and L{\'{e}}vy processes}.
\bjournal{Electron. Commun. Probab.}
\bvolume{16}
\bpages{423--434}.
\bid{mr={2831081}}
\end{barticle}
\endbibitem

\bibitem[\protect\citeauthoryear{Andersen}{1950}]{MR0039939}
\begin{barticle}[author]
\bauthor{\bsnm{Andersen},~\bfnm{Erik~Sparre}\binits{E.~S.}}
(\byear{1950}).
\btitle{On the frequency of positive partial sums of a~series of random
  variables}.
\bjournal{Mat. Tidsskr. B}
\bvolume{1950}
\bpages{33--35}.
\bid{mr={0039939}}
\end{barticle}
\endbibitem

\bibitem[\protect\citeauthoryear{Andersen}{1953a}]{MR0060173}
\begin{barticle}[author]
\bauthor{\bsnm{Andersen},~\bfnm{Erik~Sparre}\binits{E.~S.}}
(\byear{1953}a).
\btitle{On sums of symmetrically dependent random variables}.
\bjournal{Skand. Aktuarietidskr.}
\bvolume{36}
\bpages{123--138}.
\bid{mr={0060173}}
\end{barticle}
\endbibitem

\bibitem[\protect\citeauthoryear{Andersen}{1953b}]{MR0058893}
\begin{barticle}[author]
\bauthor{\bsnm{Andersen},~\bfnm{Erik~Sparre}\binits{E.~S.}}
(\byear{1953}b).
\btitle{On the fluctuations of sums of random variables}.
\bjournal{Math. Scand.}
\bvolume{1}
\bpages{263--285}.
\bid{mr={0058893}}
\end{barticle}
\endbibitem

\bibitem[\protect\citeauthoryear{Andersen}{1954}]{MR0068154}
\begin{barticle}[author]
\bauthor{\bsnm{Andersen},~\bfnm{Erik~Sparre}\binits{E.~S.}}
(\byear{1954}).
\btitle{On the fluctuations of sums of random variables. {II}}.
\bjournal{Math. Scand.}
\bvolume{2}
\bpages{195--223}.
\bid{mr={0068154}}
\end{barticle}
\endbibitem

\bibitem[\protect\citeauthoryear{Balabdaoui and
  Pitman}{2009}]{Balabdaoui2009fk}
\begin{barticle}[author]
\bauthor{\bsnm{Balabdaoui},~\bfnm{Fadoua}\binits{F.}} \AND
  \bauthor{\bsnm{Pitman},~\bfnm{Jim}\binits{J.}}
(\byear{2009}).
\btitle{The distribution of the maximal difference between {B}rownian bridge
  and its concave majorant}.
\bjournal{Bernoulli}
\bvolume{17}
\bpages{466--483}.
\bid{mr={2798000}}
\end{barticle}
\endbibitem

\bibitem[\protect\citeauthoryear{Bertoin}{1996}]{MR1406564}
\begin{bbook}[author]
\bauthor{\bsnm{Bertoin},~\bfnm{Jean}\binits{J.}}
(\byear{1996}).
\btitle{L{\'{e}}vy Processes}.
\bseries{Cambridge Tracts in Mathematics}
\bvolume{121}.
\bpublisher{Cambridge Univ. Press}, \baddress{Cambridge}.
\bid{mr={1406564}}
\end{bbook}
\endbibitem

\bibitem[\protect\citeauthoryear{Bertoin}{2000}]{MR1747095}
\begin{barticle}[author]
\bauthor{\bsnm{Bertoin},~\bfnm{Jean}\binits{J.}}
(\byear{2000}).
\btitle{The convex minorant of the {C}auchy process}.
\bjournal{Electron. Comm. Probab.}
\bvolume{5}
\bpages{51--55 (electronic)}.
\bid{mr={1747095}}
\end{barticle}
\endbibitem

\bibitem[\protect\citeauthoryear{Billingsley}{1999}]{MR1700749}
\begin{bbook}[author]
\bauthor{\bsnm{Billingsley},~\bfnm{Patrick}\binits{P.}}
(\byear{1999}).
\btitle{Convergence of Probability Measures}, \bedition{2nd} ed.
\bpublisher{Wiley}, \baddress{New York}.
\bid{mr={1700749}}
\end{bbook}
\endbibitem

\bibitem[\protect\citeauthoryear{Chaumont}{1997}]{MR1465814}
\begin{barticle}[author]
\bauthor{\bsnm{Chaumont},~\bfnm{Lo{\"{\i}}c}\binits{L.}}
(\byear{1997}).
\btitle{Excursion normalis{\'{e}}e, m{\'{e}}andre et pont pour les processus de
  {L}{\'{e}}vy stables}.
\bjournal{Bull. Sci. Math.}
\bvolume{121}
\bpages{377--403}.
\bid{mr={1465814}}
\end{barticle}
\endbibitem

\bibitem[\protect\citeauthoryear{Chaumont}{2010}]{Chaumont2010fk}
\begin{bunpublished}[author]
\bauthor{\bsnm{Chaumont},~\bfnm{Lo{\"{\i}}c}\binits{L.}}
(\byear{2010}).
\btitle{On the law of the supremum of {L}{\'{e}}vy processes}.
\bnote{Available at \url{http://arxiv.org/abs/1011.4151}}.
\end{bunpublished}
\endbibitem

\bibitem[\protect\citeauthoryear{Chaumont and Doney}{2005}]{MR2164035}
\begin{barticle}[author]
\bauthor{\bsnm{Chaumont},~\bfnm{Lo{\"{\i}}c}\binits{L.}} \AND
  \bauthor{\bsnm{Doney},~\bfnm{Ronald~A.}\binits{R.~A.}}
(\byear{2005}).
\btitle{On {L}{\'{e}}vy processes conditioned to stay positive}.
\bjournal{Electron. J. Probab.}
\bvolume{10}
\bpages{948--961 (electronic)}.
\bid{mr={2164035}}
\end{barticle}
\endbibitem

\bibitem[\protect\citeauthoryear{Chung and Fuchs}{1951}]{MR0040610}
\begin{barticle}[author]
\bauthor{\bsnm{Chung},~\bfnm{K.~L.}\binits{K.~L.}} \AND
  \bauthor{\bsnm{Fuchs},~\bfnm{W.~H.~J.}\binits{W.~H.~J.}}
(\byear{1951}).
\btitle{On the distribution of values of sums of random variables}.
\bjournal{Mem. Amer. Math. Soc.}
\bvolume{1951}
\bpages{12}.
\bid{mr={0040610}}
\end{barticle}
\endbibitem

\bibitem[\protect\citeauthoryear{Chung and Ornstein}{1962}]{MR0133148}
\begin{barticle}[author]
\bauthor{\bsnm{Chung},~\bfnm{K.~L.}\binits{K.~L.}} \AND
  \bauthor{\bsnm{Ornstein},~\bfnm{Donald}\binits{D.}}
(\byear{1962}).
\btitle{On the recurrence of sums of random variables}.
\bjournal{Bull. Amer. Math. Soc.}
\bvolume{68}
\bpages{30--32}.
\bid{mr={0133148}}
\end{barticle}
\endbibitem

\bibitem[\protect\citeauthoryear{Doney}{2007}]{MR2320889}
\begin{bbook}[author]
\bauthor{\bsnm{Doney},~\bfnm{Ronald~A.}\binits{R.~A.}}
(\byear{2007}).
\btitle{Fluctuation Theory for {L}{\'{e}}vy Processes}.
\bseries{Lecture Notes in Math.}
\bvolume{1897}.
\bpublisher{Springer}, \baddress{Berlin}.
\bid{mr={2320889}}
\end{bbook}
\endbibitem

\bibitem[\protect\citeauthoryear{Erickson}{1973}]{MR0336806}
\begin{barticle}[author]
\bauthor{\bsnm{Erickson},~\bfnm{K.~Bruce}\binits{K.~B.}}
(\byear{1973}).
\btitle{The strong law of large numbers when the mean is undefined}.
\bjournal{Trans. Amer. Math. Soc.}
\bvolume{185}
\bpages{371--381}.
\bid{mr={0336806}}
\end{barticle}
\endbibitem

\bibitem[\protect\citeauthoryear{Ethier and Kurtz}{1986}]{MR838085}
\begin{bbook}[author]
\bauthor{\bsnm{Ethier},~\bfnm{Stewart~N.}\binits{S.~N.}} \AND
  \bauthor{\bsnm{Kurtz},~\bfnm{Thomas~G.}\binits{T.~G.}}
(\byear{1986}).
\btitle{Markov Processes: Characterization and Convergence}.
\bpublisher{Wiley}, \baddress{New York}.
\bid{mr={0838085}}
\end{bbook}
\endbibitem

\bibitem[\protect\citeauthoryear{Fourati}{2005}]{MR2139029}
\begin{barticle}[author]
\bauthor{\bsnm{Fourati},~\bfnm{Sonia}\binits{S.}}
(\byear{2005}).
\btitle{Vervaat et {L}{\'{e}}vy}.
\bjournal{Ann. Inst. H. Poincar{\'{e}} Probab. Statist.}
\bvolume{41}
\bpages{461--478}.
\bid{mr={2139029}}
\end{barticle}
\endbibitem

\bibitem[\protect\citeauthoryear{Gihman and Skorohod}{1975}]{MR0375463}
\begin{bbook}[author]
\bauthor{\bsnm{Gihman},~\bfnm{{\u{I}}.~{\={I}}.}\binits{{\u{I}}.~{\={I}}.}} \AND \bauthor{\bsnm{Skorohod},~\bfnm{A.~V.}\binits{A.~V.}}
(\byear{1975}).
\btitle{The Theory of Stochastic Processes. {II}}.
\bpublisher{Springer}, \baddress{New York}.
\bid{mr={0375463}}
\end{bbook}
\endbibitem

\bibitem[\protect\citeauthoryear{Greenwood and Pitman}{1980}]{MR588409}
\begin{barticle}[author]
\bauthor{\bsnm{Greenwood},~\bfnm{Priscilla}\binits{P.}} \AND
  \bauthor{\bsnm{Pitman},~\bfnm{Jim}\binits{J.}}
(\byear{1980}).
\btitle{Fluctuation identities for {L}{\'{e}}vy processes and splitting at the
  maximum}.
\bjournal{Adv. in Appl. Probab.}
\bvolume{12}
\bpages{893--902}.
\bid{mr={0588409}}
\end{barticle}
\endbibitem

\bibitem[\protect\citeauthoryear{Groeneboom}{1983}]{MR714964}
\begin{barticle}[author]
\bauthor{\bsnm{Groeneboom},~\bfnm{Piet}\binits{P.}}
(\byear{1983}).
\btitle{The concave majorant of {B}rownian motion}.
\bjournal{Ann. Probab.}
\bvolume{11}
\bpages{1016--1027}.
\bid{mr={0714964}}
\end{barticle}
\endbibitem

\bibitem[\protect\citeauthoryear{Kallenberg}{2002}]{MR1876169}
\begin{bbook}[author]
\bauthor{\bsnm{Kallenberg},~\bfnm{Olav}\binits{O.}}
(\byear{2002}).
\btitle{Foundations of Modern Probability}, \bedition{2nd} ed.
\bpublisher{Springer}, \baddress{New York}.%
\bid{mr={1876169}}
\end{bbook}
\endbibitem

\bibitem[\protect\citeauthoryear{Knight}{1996}]{MR1417982}
\begin{barticle}[author]
\bauthor{\bsnm{Knight},~\bfnm{F.~B.}\binits{F.~B.}}
(\byear{1996}).
\btitle{The uniform law for exchangeable and {L}{\'{e}}vy process bridges}.
\bjournal{Ast{\'{e}}risque}
\bvolume{236}
\bpages{171--188}.
\bid{mr={1417982}}
\end{barticle}
\endbibitem

\bibitem[\protect\citeauthoryear{Kyprianou}{2006}]{110460001}
\begin{bbook}[author]
\bauthor{\bsnm{Kyprianou},~\bfnm{Andreas~E.}\binits{A.~E.}}
(\byear{2006}).
\btitle{{Introductory Lectures on Fluctuations of L\'{e}vy Processes with
  Applications.}}
\bpublisher{Springer}, \baddress{Berlin}.
\end{bbook}
\endbibitem

\bibitem[\protect\citeauthoryear{Lachieze-Rey}{2009}]{Lachieze-Rey2009fk}
\begin{bmisc}[author]
\bauthor{\bsnm{Lachieze-Rey},~\bfnm{Raphael}\binits{R.}}
(\byear{2009}).
\bhowpublished{Concave majorant of stochastic processes and {B}urgers turbulence.
\textit{J.~Theor. Probab.} To appear.
DOI:\doiurl{10.1007/S10959-011-0354-7}.}
\end{bmisc}
\endbibitem

\bibitem[\protect\citeauthoryear{McCloskey}{1965}]{McCloskey}
\begin{bmisc}[author]
\bauthor{\bsnm{McCloskey},~\bfnm{J.~W.}\binits{J.~W.}}
(\byear{1965}).
\bhowpublished{A model for the distribution of individuals by species in an
  environment.
Ph.D. thesis, Michigan State Univ., Ann Anbor, MI.}
\bid{mr={2615013}}
\end{bmisc}
\endbibitem

\bibitem[\protect\citeauthoryear{Miermont}{2001}]{MR1844511}
\begin{barticle}[author]
\bauthor{\bsnm{Miermont},~\bfnm{Gr{{\'{e}}}gory}\binits{G.}}
(\byear{2001}).
\btitle{Ordered additive coalescent and fragmentations associated to {L}evy
  processes with no positive jumps}.
\bjournal{Electron. J. Probab.}
\bvolume{6}
\bpages{33 pp. (electronic)}.
\bid{mr={1844511}}
\end{barticle}
\endbibitem

\bibitem[\protect\citeauthoryear{Millar}{1977}]{MR0433606}
\begin{barticle}[author]
\bauthor{\bsnm{Millar},~\bfnm{P.~W.}\binits{P.~W.}}
(\byear{1977}).
\btitle{Zero--one laws and the minimum of a~{M}arkov process}.
\bjournal{Trans. Amer. Math. Soc.}
\bvolume{226}
\bpages{365--391}.
\bid{mr={0433606}}
\end{barticle}
\endbibitem

\bibitem[\protect\citeauthoryear{Nagasawa}{2000}]{MR1739699}
\begin{bbook}[author]
\bauthor{\bsnm{Nagasawa},~\bfnm{Masao}\binits{M.}}
(\byear{2000}).
\btitle{Stochastic Processes in Quantum Physics}.
\bseries{Monographs in Mathematics}
\bvolume{94}
\bpages{355--388}.
\bpublisher{Birkh{\"{a}}user}, \baddress{Basel}.
\bid{mr={1739699}}
\end{bbook}
\endbibitem

\bibitem[\protect\citeauthoryear{Pe{\v{c}}erski{\u\i} and
  Rogozin}{1969}]{MR0260005}
\begin{barticle}[author]
\bauthor{\bsnm{Pe{\v{c}}erski{\u\i}},~\bfnm{E.~A.}\binits{E.~A.}} \AND
  \bauthor{\bsnm{Rogozin},~\bfnm{B.~A.}\binits{B.~A.}}
(\byear{1969}).
\btitle{The combined distributions of the random variables connected with the
  fluctuations of a~process with independent increments}.
\bjournal{Teor. Veroyatn. Primen.}
\bvolume{14}
\bpages{431--444}.
\bid{mr={0260005}}
\end{barticle}
\endbibitem

\bibitem[\protect\citeauthoryear{Perman, Pitman and Yor}{1992}]{MR1156448}
\begin{barticle}[author]
\bauthor{\bsnm{Perman},~\bfnm{Mihael}\binits{M.}},
  \bauthor{\bsnm{Pitman},~\bfnm{Jim}\binits{J.}} \AND
  \bauthor{\bsnm{Yor},~\bfnm{Marc}\binits{M.}}
(\byear{1992}).
\btitle{Size-biased sampling of {P}oisson point processes and excursions}.
\bjournal{Probab. Theory Related Fields}
\bvolume{92}
\bpages{21--39}.
\bid{mr={1156448}}
\end{barticle}
\endbibitem

\bibitem[\protect\citeauthoryear{Pitman}{1983}]{MR733673}
\begin{bincollection}[author]
\bauthor{\bsnm{Pitman},~\bfnm{J.~W.}\binits{J.~W.}}
(\byear{1983}).
\btitle{Remarks on the convex minorant of {B}rownian motion}.
In \bbooktitle{Seminar on Stochastic Processes, 1982 ({E}vanston, IL,
  1982)}.
\bseries{Progress in Probability Statist.}
\bvolume{5}
\bpages{219--227}.
\bpublisher{Birkh{\"{a}}user}, \baddress{Boston, MA}.
\bid{mr={0733673}}
\end{bincollection}
\endbibitem

\bibitem[\protect\citeauthoryear{Pitman and Ross}{2010}]{2010PitmanRoss}
\begin{bmisc}[author]
\bauthor{\bsnm{Pitman},~\bfnm{Jim}\binits{J.}} \AND
  \bauthor{\bsnm{Ross},~\bfnm{Nathan}\binits{N.}}
(\byear{2010}).
\bhowpublished{The greatest convex minorant of Brownian motion, meander, and bridge.
\textit{Probab. Theory Related Fields.} To appear.
DOI:\href{http://dx.doi.org/10.1007/S00440-011-0385-0}{10.1007/}
\href{http://dx.doi.org/10.1007/S00440-011-0385-0}{S00440-011-0385-0}}.
\end{bmisc}
\endbibitem

\bibitem[\protect\citeauthoryear{Pitman and Yor}{1997}]{MR1434129}
\begin{barticle}[author]
\bauthor{\bsnm{Pitman},~\bfnm{Jim}\binits{J.}} \AND
  \bauthor{\bsnm{Yor},~\bfnm{Marc}\binits{M.}}
(\byear{1997}).
\btitle{The two-parameter {P}oisson--{D}irichlet distribution derived from a~stable subordinator}.
\bjournal{Ann. Probab.}
\bvolume{25}
\bpages{855--900}.
\bid{mr={1434129}}
\end{barticle}
\endbibitem

\bibitem[\protect\citeauthoryear{Rogozin}{1968}]{MR0242261}
\begin{barticle}[author]
\bauthor{\bsnm{Rogozin},~\bfnm{B.~A.}\binits{B.~A.}}
(\byear{1968}).
\btitle{The local behavior of processes with independent increments}.
\bjournal{Teor. Veroyatn. Primen.}
\bvolume{13}
\bpages{507--512}.
\bid{mr={0242261}}
\end{barticle}
\endbibitem

\bibitem[\protect\citeauthoryear{Sato}{1999}]{MR1739520}
\begin{bbook}[author]
\bauthor{\bsnm{Sato},~\bfnm{Ken-iti}\binits{K.-i.}}
(\byear{1999}).
\btitle{L{\'{e}}vy Processes and Infinitely Divisible Distributions}.
\bseries{Cambridge Stud. Adv. Math.}
\bvolume{68}.
\bpublisher{Cambridge Univ. Press}, \baddress{Cambridge}.
\bid{mr={1739520}}
\end{bbook}
\endbibitem

\bibitem[\protect\citeauthoryear{Spitzer}{1956}]{MR0079851}
\begin{barticle}[author]
\bauthor{\bsnm{Spitzer},~\bfnm{Frank}\binits{F.}}
(\byear{1956}).
\btitle{A combinatorial lemma and its application to probability theory}.
\bjournal{Trans. Amer. Math. Soc.}
\bvolume{82}
\bpages{323--339}.
\bid{mr={0079851}}
\end{barticle}
\endbibitem

\bibitem[\protect\citeauthoryear{Suidan}{2001a}]{MR1978665}
\begin{barticle}[author]
\bauthor{\bsnm{Suidan},~\bfnm{T.~M.}\binits{T.~M.}}
(\byear{2001}a).
\btitle{Convex minorants of random walks and {B}rownian motion}.
\bjournal{Teor. Veroyatn. Primen.}
\bvolume{46}
\bpages{498--512}.
\bid{mr={1978665}}
\end{barticle}
\endbibitem

\bibitem[\protect\citeauthoryear{Suidan}{2001b}]{MR1861441}
\begin{barticle}[author]
\bauthor{\bsnm{Suidan},~\bfnm{T.~M.}\binits{T.~M.}}
(\byear{2001}b).
\btitle{A one-dimensional gravitationally interacting gas and the convex
  minorant of {B}rownian motion}.
\bjournal{Uspekhi Mat. Nauk}
\bvolume{56}
\bpages{73--96}.
\bid{mr={1861441}}
\end{barticle}
\endbibitem

\bibitem[\protect\citeauthoryear{Uribe~Bravo}{2011}]{UribeBravoVervaatExtended}
\begin{bmisc}[author]
\bauthor{\bsnm{Uribe~Bravo},~\bfnm{G.}\binits{G.}}
(\byear{2011}).
\btitle{Bridges of {L}\'{e}vy processes conditioned to stay positive.}
\bhowpublished{Available at \url{http://arxiv.org/abs/1101.4184}}.
\end{bmisc}
\endbibitem

\bibitem[\protect\citeauthoryear{Vervaat}{1979}]{MR515820}
\begin{barticle}[author]
\bauthor{\bsnm{Vervaat},~\bfnm{Wim}\binits{W.}}
(\byear{1979}).
\btitle{A relation between {B}rownian bridge and {B}rownian excursion}.
\bjournal{Ann. Probab.}
\bvolume{7}
\bpages{143--149}.
\bid{mr={0515820}}
\end{barticle}
\endbibitem

\bibitem[\protect\citeauthoryear{Vigon}{2002}]{MR1875147}
\begin{barticle}[author]
\bauthor{\bsnm{Vigon},~\bfnm{Vincent}\binits{V.}}
(\byear{2002}).
\btitle{Votre {L}{\'{e}}vy rampe-t-il?}
\bjournal{J. Lond. Math. Soc. (2)}
\bvolume{65}
\bpages{243--256}.
\bid{mr={1875147}}
\end{barticle}
\endbibitem

\end{thebibliography}
%

%


\printaddresses

\end{document}